\documentclass{amsart}
\usepackage{amsmath, amsthm, amsopn, amsfonts}
\usepackage{graphicx}
\usepackage{graphpap}

\def\pasdegrille{\let\grille = \pasgrille}
\def\ecriture#1#2{\setbox1=\hbox{#1}
\dimen1= \wd1
\dimen2=\ht1
\dimen3=\dp1
\grille #2 \box1 }
\def\aat#1#2#3{
\divide \dimen1 by 48
\dimen3=\dimen1
\multiply \dimen1 by #1
\advance \dimen1 by -\dimen3
\divide \dimen1 by 101
\multiply \dimen1 by 100
\divide \dimen2 by \count11
\multiply \dimen2 by #2 
\setbox0=\hbox{#3}\ht0=0pt\dp0=0pt
  \rlap{\kern\dimen1 \vbox to0pt{\kern-\dimen2\box0\vss}}\dimen1= \wd1
\dimen2=\ht1}
\def\pasgrille{
\count12= \dimen1 
\divide \count12 by 50
\divide \dimen2 by \count12
\count11 =\dimen2
\ 
\divide \dimen1 by 48
\setlength{\unitlength}{\dimen1}
\smash{\rlap{\ }}
\dimen1= \wd1
\dimen2=\ht1
}
\def\grille{
\count12= \dimen1 
\divide \count12 by 50
\divide \dimen2 by \count12
\count11 =\dimen2
\ 
\divide \dimen1 by 48
\setlength{\unitlength}{\dimen1}
\smash{\rlap{\graphpaper[1](0,0)(50, \count11)}}
\dimen1= \wd1
\dimen2=\ht1
}

\pasdegrille
\usepackage{amssymb}
\ifx\pdfoutput\undefined
\usepackage{graphicx}  \DeclareGraphicsExtensions{.ps}

\else
\usepackage{amsfonts}
\usepackage[pdftex]{graphicx}  \DeclareGraphicsExtensions{.pdf,.mps}

\fi

\usepackage[T1]{fontenc}

\usepackage{epic, eepic}
\def\squarebox#1{\hbox to #1{\hfill\vbox to #1{\vfill}}} 
\newcommand{\stopthm}{\hfill\hfill\vbox{\hrule\hbox{\vrule\squarebox 
                 {.667em}\vrule}\hrule}\smallskip} 

\newcommand{\1}{{\bold 1}}

\newcommand{\CI}{{\mathcal C}^\infty }
\newcommand{\HH}{{\mathcal H}}
\newcommand{\CIc}{{\mathcal C}^\infty_{\rm{c}} }
\newcommand{\CIb}{{\mathcal C}^\infty_{\rm{b}} }
\newcommand{\Op}{{\operatorname{Op}^{{w}}_h}}

\newcommand{\RR}{{\mathbb R}}
\newcommand{\SP}{{\mathbb S}}

\newcommand{\NN}{{\mathbb N}}

\newcommand{\supp}{\operatorname{supp}}

\newcommand{\rest}{\!\!\restriction}

\renewcommand{\Re}{\mathop{\rm Re}\nolimits}
\renewcommand{\Im}{\mathop{\rm Im}\nolimits}

\theoremstyle{plain}

\newtheorem{thm}{Theorem}
\newtheorem{prop}{Proposition}[section]
\newtheorem{lem}{Lemma}[section]

\theoremstyle{definition}

\newtheorem{rem}{Remark}

\numberwithin{equation}{section}

\def\bbbone{{\mathchoice {1\mskip-4mu \rm{l}} {1\mskip-4mu \rm{l}}
{ 1\mskip-4.5mu \rm{l}} { 1\mskip-5mu \rm{l}}}}

\title[Control in the presence of a black box]
{Geometric control in the presence of a black box}
\author[N. Burq]{Nicolas Burq}
\address{Universit{\'e} Paris Sud,
Math{\'e}matiques,
B{\^a}t 425, 91405 Orsay Cedex, France}
\email{Nicolas.burq@math.u-psud.fr}
\author[M. Zworski]{Maciej Zworski}
\address{Mathematics Department, University of California \\
Evans Hall, Berkeley, CA 94720, USA}
\email{zworski@math.berkeley.edu}

\usepackage{amssymb}
\usepackage{amsmath, amsthm, amsopn, amsfonts}

\def\11{{\rm 1~\hspace{-1.4ex}l} }

\begin{document}

\begin{abstract}
We apply the  ``black box'' scattering theory
to problems in control theory for the Schr\"odinger equation,
and in high energy eigenvalue scarring.
\end{abstract}
   
\maketitle   
   
\section{Introduction}   
\label{in}

The purpose of this paper
is to show how ideas coming from scattering theory (resolvent estimates) 
lead to  results in control theory and to some closely related eigenfunction 
estimates. 

The black box approach in scattering theory developed by 
Sj\"ostrand and the second author \cite{SjZw91} 
puts scattering problems with different structures
in one framework, and allows abstract applications of spectral 
results known for confined systems. One striking example is 
a reduction of scattering on finite volume surfaces to one dimensional
black box scattering. In this paper we take the opposite point
of view: a black box in a confined system is replaced by 
a scattering problem. That permits having isolated dynamical 
phenomena (such as only one closed orbit) impossible in confined
systems. It also permits using some finer results of scattering theory
directly. 

We stress that this follows the well established trend 
(see Bardos-Lebeau-Rauch \cite{BLR}) of using propagation 
of singularities results developed for scattering theory in geometric 
control theory. We also mention that the term ``black box'' is
commonly used, in a similar context, in applied control theory 
\cite{BB}.

Since the proofs are simple and since it is profitable
to state the results in an abstract setting which requires 
a certain amount of preparation, in this section we will 
present some typical applications.

\begin{figure}[ht]
$$\ecriture{\includegraphics[width=12cm]{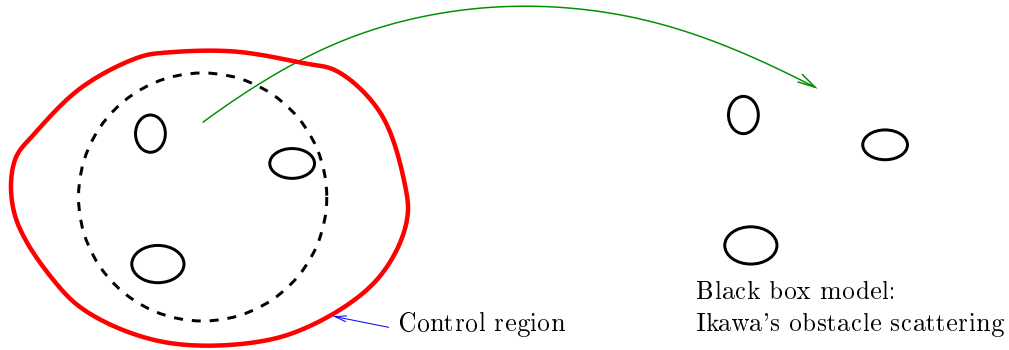}}
{\aat{22}{1}{Control region}\aat{38}{1}{\vbox{\hbox{Black box model:} \hbox{Ikawa's obstacle scattering}}}}
$$
\caption{Control in the exterior of several convex bodies}\label{fig:burq}
\end{figure}

In geometric control theory for the Schr\"odinger equation 
(see Lebeau \cite{Le92}, and also \cite{Li88},\cite{Zu88} 
for earlier work and background) we are concerned with the following mixed 
problem:
\begin{gather}
\label{eq:1.1}\begin{gathered}
 (i\partial_{t}+ \Delta)u=0 \ \ \text{ in $ \Omega $ }\\
\begin{split}
&  u\rest_{[0,T] \times \partial \Omega }= g\bbbone_{[0,T]\times \Gamma}\\
& u\rest_{t=0}=u_{0}.
\end{split}
\end{gathered}
\end{gather}where $\Omega$ is an open subset of ${\mathbb R}^d$, $\partial \Omega$ is its boundary and $\Gamma$ is an open subset of $\partial \Omega$. 
The question is to determine a (large) class of functions $u_{0}$ for which there exists
a {\em control $g$} such that $u\rest_{t=T}=0$.
In a geometric setting in which full geometric control 
fails, the following result was established by the first author in 
\cite{Bu93}:
\begin{thm}
\label{t:burq}
Consider $\Theta= \cup_{j=1}^N \Theta_{j} \subset \RR^d $, a union of mutually 
disjoint closed sets with 
strictly convex smooth boundaries, and  satisfying the assumptions in Sect.\ref{ikba} below. 
Let $\widetilde { \Omega}$ be a bounded domain with a smooth boundary 
and containing ${\rm{convhull}}(\Theta)$. 
Denote by $\Omega= \widetilde \Omega \setminus \Theta$ and 
$\Gamma= \partial \widetilde{ \Omega}$. 
Then for any $T, \varepsilon>0$ and any $u_{0}\in H^{1+\varepsilon}_{0}( \Omega)$ there 
exists $g\in L^2( [0,T] \times \Gamma)$ such that in \eqref{eq:1.1}, with we have  
$ u \rest_{t>T}\equiv 0$. 
\end{thm}
In Fig.\ref{fig:burq} on the left we have three convex obstacles inside of the
{boundary} of $ \widetilde \Omega $. Inside of the black box 
bounded by the dotted line the local geometry is the same as in the scattering
problem on the right.

We are going to show how 
Theorem \ref{t:burq} can be obtained directly from 
estimates on the resolvent of the Laplace operator, which in turn can be deduced 
from semi-classical microlocal analysis or from known results in scattering 
theory. In the case quoted above, these come from the work of Ikawa \cite{Ik}
and in particular we can now avoid most of the delicate analysis of~\cite{Bu93}.

The next application generalizes a result of Colin de Verdi\`ere and Parisse \cite{CdVP} 
who considered a special case of an isolated trajectory lying on a segment of
a constant negative curvature cylinder in dimension two:
\begin{thm}
\label{t:cdvp}
Suppose that $ ( X , g ) $ is a compact Riemannian manifold with a (possibly empty) 
boundary and $ \gamma \subset X $ is  
a closed hyperbolic geodesic (we allow broken geodesic flow as long as 
the reflections are all transversal).
If $ \chi \in \CI ( X , [0, 1] ) $ is supported
in a sufficiently small neighbourhood of $ \gamma $ then there exists a constant
$ C = C ( \gamma ) $ such that for {\em any} 
eigenfunction, $ u $, of the Laplacian, $ \Delta_g $ with Dirichlet or Neumann boundary conditions, 
we have 
\begin{equation}
\label{eq:1.cdv}
C \int_X | u (x) |^2 ( 1 - \chi) ( x ) d {\rm{vol}}_g \geq \frac{1}{\log \lambda }
\int_X | u (x) |^2  d {\rm{vol}}_g \,, \ \ -\Delta_g u = \lambda u 
\,.
\end{equation}
\end{thm}
An example \cite{CdVP} of a cylinder segment with Dirichlet boundary conditions
shows that the result is optimal. 

The proof of Theorem \ref{t:cdvp} (see also Theorem \ref{t:cdvp}$'$) 
is based on putting the closed hyperbolic orbit into a
{\em microlocal black box}, where that orbit becomes the only trapped orbit in a
scattering problem. We can then use scattering estimates based on the quantum 
monodromy method \cite{SjZw02}, and the work of 
G\'erard \cite{Ge} and G\'erard-Sj\"ostrand \cite{GeSj}
to obtain estimates leading to 
\eqref{eq:1.cdv}.

\begin{figure}[ht]
\includegraphics[width=10cm]{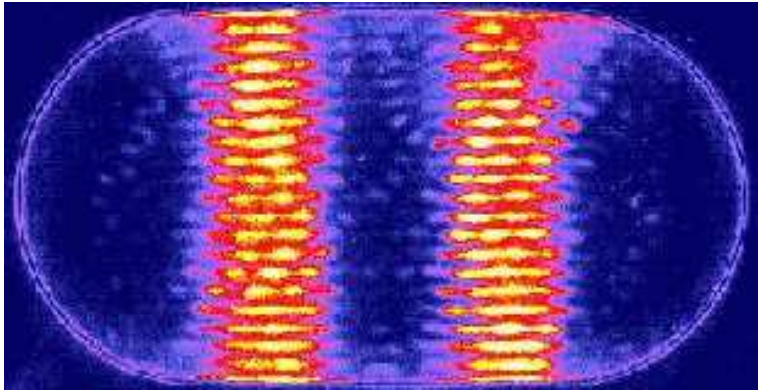}
\caption{An experimental image of the wave in the 
``black box'' in Fig. \ref{fig:bun} -- 
see \cite{CH96} and 
{\tt http://www.bath.ac.uk/$\sim$pyscmd/acoustics}.}
\label{fig:bath}
\end{figure}

We conclude with a brief discussion of another example related to {\em
eigenvalue scarring} (see Theorem \ref{t:bun} below for a full discussion).
While in Theorem \ref{t:cdvp} we eliminated the need for separation of 
variables, its use is essential in this case.
For the Bunimovich cavity shown in Fig.\ref{fig:bath} the natural
black box for constructing bouncing ball modes
(two are shown in the same figure)
is a rectangle constituting the central part of the cavity -- see the
recent discussion of this in \cite{Do} and \cite{Ze}. On one hand, our result shows that
the crude error estimate
\begin{equation}
\label{eq:1.bun}
   ( - \Delta_D - \lambda ) u_\lambda = {\mathcal O}(1) \,, \ \ \|u \|=1 \,, 
\end{equation}
in the quasimodes obtained by truncating the rectangle modes is in fact
the best possible and on the other hand that the eigenfunctions can {\em not} acumulate at high frequency only in the central part. This agrees with 
the experimental results \cite{CH96} where it was stressed that
phenomena shown in  Fig.\ref{fig:bath} can occur only at low frequencies
(see also \cite{BS} for a different discussion and references to the physics
literature). For 
an exact eigenstate we have the following 
\begin{thm}
\label{t:bunim}
Let $u$ be a Dirichlet eigenfunction of the Laplacian on 
the Bunimovich stadium $M$: 
\begin{equation*}-\Delta u= \lambda u , \qquad u\rest_{\partial M}=0
\end{equation*} 
Let $ a(x) $ be any continuous function 
identically $ 1 $ on the non-rectangular part of $M$. Then there exists $ C > 0 $ such that
\begin{equation}\label{eq:1.bunimo}
C \int_{M }
| a(x) u (x)|^2 dx \geq \int_M |u ( x) |^2 dx \,.
\end{equation}
\end{thm}
Stronger results (implying \eqref{eq:1.bunimo})
are presented in Theorems \ref{t:bunim}$'$ and 
\ref{t:bun} in Sect.\ref{s:buni}. A self contained proof 
of Theorem \ref{t:bunim} and a discussion of related mathematical
and physical literature has been presented in \cite{BZ3}.
We stress that only the properties of the rectangular part used 
as a ``black box'' are needed for this result. 

\medskip
\noindent
{\sc Acknowledgments.} The authors would like to thank the National 
Science Foundation for partial support under the grant  DMS-0200732.
They are also grateful to Steve Zelditch for informing them of \cite{Do} and
\cite{Ze} which expanded the breadth of this note, to Luc Miller for
helpful comments on the first version of the paper, and to 
Victor Humphrey and Paul Chinnery for the permission to use their Fig.\ref{fig:bath}. 
The first author thanks the Mathematical Science Research Institute for its hospitality 
during spring 2003.

\section{Preliminaries}
\label{pr} 
In this section we review some basic aspects of semiclassical 
microlocal analysis, following \cite[Section 3]{SjZw02}. 
Thus, let $ X $ be a compact $\CI$ manifold. 
We consider pseudo-differential operators as acting
on half-densities, 
$  u ( x ) |dx|^{\frac12} \in \CI ( X , \Omega_X^{\frac12} ) $,
where we use the informal notation indicating how the half-densities 
change under changes of variables: 
\[ u ( x ) |d x |^{\frac12} =  v ( y ) |d y|^{\frac12}\,,
\  y = \kappa ( x ) \ \Longleftrightarrow \
v  ( \kappa ( x ) ) | \kappa ' ( x ) |^{\frac12 } = u ( x  ) \,,
\]
Consequently the symbols will also be considered
as half-densities -- see \cite[Sect.18.1]{Hor} for a general introduction
and \cite[Appendix]{SjZw02} for a  discussion of the semi-classical case.
This way our results are more general and do not depend on the
choice of a metric on $ X$. If $ X $ is a Riemannian manifold and the
operator we consider its Laplace-Bertrami operator then the natural
Riemannian density is all we need.

By
symbols on $ X $ we mean the following class:
\[ S^{k,m} ( T^* X , \Omega_{T^*X}^{\frac12}
) = \{ a \in \CI( T^* X \times (0, 1] ; \Omega_{T^*X}^{\frac12} ) :
|\partial_x ^{ \alpha } \partial _\xi^\beta a ( x, \xi ;h ) | \leq
C_{ \alpha, \beta} h^{ -m } \langle \xi \rangle^{ k - |\beta| } 
\} \,, \]
and the class  corresponding pseudodifferential operators, $ 
\Psi_{h}^{m,k} ( X , \Omega_X^{\frac12}) $,
obtained from a local formula in $ \RR^n $:
\begin{equation}
\label{eq:weyl}
 \Op (a) u ( x) = \frac{1}{ ( 2 \pi h )^n } 
\int \int  a \left( \frac{x + y }{2}  , \xi,h \right
) e^{ i \langle x -  y, \xi \rangle / h } u ( y ) 
dy d \xi \,. \end{equation}
The principal symbol map, 
\[  \sigma_h \; :\; \Psi_{h}^{m,k} ( X , \Omega_X^{\frac12}) 
\; \longrightarrow \; S^{k,m} / S^{k-1,m-1} ( T^* X , \Omega_{T^*X}^{\frac12})\,, \]
gives the left inverse of $ \Op $ in the sense that 
$ \sigma_h \circ \Op : S^{m,k} \rightarrow 
S^{m,k}/S^{m-1,k-1} $ is the natural projection. 
We refer to \cite{DiSj} for a detailed discussion of the Weyl quantization
and to \cite{WZ} for a discussion in the case of manifolds.

For $ a \in S^{m,k} ( T^* X \Omega_{T^*X}^{\frac12}) $ we follow \cite{SjZw02} in defining 
\[ {\text{ess-supp}}_h\; a \subset T^* X \sqcup S^* X \,, \ \ 
S^* X \stackrel{{\rm{def}}}{=}  ( T^* X  \setminus 0 ) / {\mathbb R}_+ \,, 
\]
where the usual $ \RR_+ $ action is given by multiplication on 
the fibers: $ ( x, \xi ) \mapsto  ( x , t \xi ) $, as 
\[
\begin{split}
&  {\text{ess-supp}}_h\; a =  \\
& \ \ \ \complement \{ ( x, \xi ) \in T^* X \; : \; 
\exists \; \epsilon > 0 \ \partial_x ^\alpha \partial_\xi ^\beta 
a ( x', \xi' ) = {\mathcal O} ( h^\infty ) 
  \,, \  d( x, x' ) + | \xi - \xi' |  < \epsilon \} \\
& \ \ \ \cup \complement  \{ 
( x, \xi ) \in T^* X \setminus 0 \; : \; 
\exists \; \epsilon > 0 \ \partial_x ^\alpha \partial_\xi ^\beta 
a ( x', \xi' ) = {\mathcal O} ( h^\infty  \langle \xi' \rangle^{  -\infty}) 
  \,,  \\ 
& \ \ \ \ \ \ \ \ \ 
d( x, x' ) + 1 / |\xi'| + | \xi/ |\xi|  - \xi'/ |\xi'| |  
< \epsilon \} / {\mathbb R}_+ 
\end{split} \,. \]
For $ A \in \Psi^{m,k} _h ( X, \Omega_X^{\frac12} ) $, then put
\[ WF_h ( A) =  {\text{ess-supp}}_h\; a \,, \ \ A = \Op ( a ) \,, 
\]
and this definition does not depend on the choice of
$ \Op $. For 
\[  u \in \CI ( ( 0 , 1]_h ; {\mathcal D}' ( X, \Omega_X^{\frac12}) ) \,, \ \  \exists \; 
N_0\,,  \ \ h^{- N_0} u \ \text{ is bounded in $ {\mathcal D}' ( X, 
\Omega_X^{\frac12} ) $,}\]
 we 
define the {\em semi-classical wave front set} as
\[ WF_
h ( u ) = \complement \{ ( x, \xi ) \; : \; 
\exists \;  A \in \Psi_h^{0,0} ( X, \Omega_X^{\frac12}
) \ \sigma_h ( A ) ( x, \xi ) \neq 0 \,,
\ A u \in h^\infty  \CI ( ( 0 , 1]_h ; \CI ( X, \Omega_X^{\frac12}) ) \} \,.\]
When $ u $ is not necessarily smooth we can give a definition 
analogous to that of $ {\text{ess-supp}}_h\; a $.
In this paper we will work in a pure semi-classical 
setting and consequently
only {\em compact} subsets of $ T^* X $ will be important. 
Consequently, this definition is sufficient for our purposes.

We also need to review the notion of microlocal equivalence of
operators and other objects. Suppose that
\[ T  \; : \; \CI ( X, \Omega_X^{\frac12} ) \rightarrow \CI ( X, \Omega_X^{\frac12} ) \,, \]
and that for any semi-norm $ \| \bullet\|_1 $ on $ \CI ( X, \Omega_X^{\frac12}) $ there exists a semi-norm $ \| \bullet\|_2 $ and $ M_0$ such that
\[  \| T u \|_1  = {\mathcal O } ( h ^{ - M_0 } ) 
\| u \|_2 \,. \]
This condition makes $ T $ {\em semi-classically tempered}. In the sequel all operators considered will be assumed to satisfy this temperence condition.
For open sets, $ V \subset T^* X$, $ U \subset T^* X $, the operators 
{\em defined microlocally} near  $ V \times U $ 
are given by equivalence classes of tempered operators
given by the relation
\[ T \sim T' \ \Longleftrightarrow \ A ( T - T' ) B = {\mathcal O} ( h^\infty )
\; : \; {\mathcal D}' ( X, \Omega_X^{\frac12} ) \ \longrightarrow \CI ( X ,
\Omega_X^{\frac12}  ) \,, \]
for any 
$ A, B \in \Psi_h ^{0,0 } ( X , \Omega_X^{\frac12}) $ such that 
\begin{gather}
\label{eq:2.7}
\begin{gathered}
WF ( A ) \subset \widetilde V \,, \ \ WF ( B ) \subset \widetilde U\,,
\\
\bar  V \Subset  \widetilde V \Subset T^* X \,, \ \ 
\bar U \Subset  \widetilde U \Subset T^* X \,, \ \ \widetilde U \,,\,
\widetilde V \ \text{ \  open}\,.
\end{gathered}
\end{gather}

We say that $ P = Q $ {\em microlocally} near $ U \times V $ if 
$  A P B - A Q B = {\mathcal O}_{L^2 \rightarrow L^2 }
 ( h ^ \infty ) $, where because of the assumed pre-compactness of 
$ U $ and $ V$ the $ L^2 $ norms can be replaced by any other norms.
For operator identities this will 
be the meaning of equality of operators in this paper, with $ U , V $
specified (or clear from the context). Similarly, we say that $ 
B = T^{-1} $ microlocally near $ V \times V $, if $ B T = I$ microlocally 
near $ U \times U $, and $ T B = I $ microlocally near $ V \times U $.
More generally, we could say that $ P = Q $ microlocally on $ 
W \subset T^* X \times T^* X $ (or, say, $ P $ is 
microlocally defined there), if for any $ U, V$, $ U \times V \subset
W $, $ P = Q $ microlocally in $ U \times V$. 
We should stress that ``microlocally'' is always meant in this 
semi-classical sense in our paper.

Rather than review the definition of $h$-Fourier integral operators
we will recall a characterization which is essentially a converse
of Egorov's theorem:
\begin{prop}
\label{p:2.0}
Suppose that $ U = {\mathcal O} ( 1) : L^2 ( X) \rightarrow 
L^2 ( X ) $, and that for every $ A \in \Psi_h ^{0,0} ( X ) $ we have
\[ A U = U B \,, \ \ B \in \Psi^{0,0}_h ( X) \,, \ \ 
\sigma ( B ) = {\kappa}^* \sigma ( A ) \,,\]
microlocally near $ ( m_0 , m_0 )$ 
where $ \kappa : T^* X \rightarrow T^*X $ is a symplectomorphism,
defined locally near $ m_0 $, $ \kappa ( m_0 ) = m_0 $. 
Then $ U $ is, microlocally, near 
$ (m_0, m_0 )$, an $h$-Fourier integral operator of order zero, quantizing 
 $ \kappa $, that is associated to the graph of $ \kappa $.
\end{prop}
For the proof  and further details we refer the reader
to \cite[Lemma 3.4]{SjZw02}.
We will use the following well known fact (see \cite[Proposition 3.5]{SjZw02}
for the proof):
\begin{prop}
\label{p:2.1}
Suppose that $ P \in \Psi_{h}^{0,k} ( X ) $ has a real principal
symbol which 
satisfies the condition
\[ p = 0 \Longrightarrow dp \neq 0   \,.\]
For any $ m_0 \in p^{-1}(0) $
 there exists an $h$-Fourier Integral Operator,  $ F $,
\begin{gather*}
F P = hD_{x_1}  F \,, \ \ \text{microlocally near $( ( 0 , 0 ),  m_0  )$}\, \\
F^{-1} \ \ \text{exists microlocally near $(m_0 , ( 0 , 0 ) ) $}\,.
\end{gather*}
\end{prop}

\section{From resolvent estimates to time dependent control}
\label{res}
In this section we will present a simple abstract argument 
showing how semi-classical resolvent estimates give a
control result for the semi-classical Schr\"odinger operator.
An adaptation of this argument to the classical control 
setting will be presented in Sect.\ref{cbb}.

\begin{thm}
\label{t:3}
Let $ P ( h ) $ be a family of self-adjoint operators on a
Hilbert space $ \HH$, with a fixed domain $ {\mathcal D} $.
Let $ {\mathcal H}_1 $ be 
another Hilbert space, and
suppose that for a bounded family operators, $ A ( h ) : 
{\mathcal D} \rightarrow  \HH_1 $, we have
\begin{equation}
\label{eq:3.1}
  \| u \|_\HH \leq \frac{G( h ) }{h} \| (P(h) + \tau) u \|_\HH + g ( h ) \| 
A ( h ) u \|_{\HH_1}  \,, 
\end{equation}
$ \tau \in I = ( -b , -a  ) \Subset \RR $, $1\leq  G ( h ) = 
 {\mathcal O} ( h^{-N_0} )  $, for some $ N_0$. Fix $ \chi \in \CIc ((a,b) ) $.
There exists constants $ c_0 $, and $ C_0 $, such that
for any $ T ( h ) $ satisfying
\begin{equation}
\label{eq:3.2}
 \frac{ G ( h )} { T ( h )} < c_0
\end{equation}
we have for $ 0 < h < h_0( \delta ) $, 
\begin{equation}
\label{eq:3.3}
\| \chi ( P (h) ) u \|^2_\HH \leq C_0 \frac{ g( h )^2 }{ T( h ) }
\int_0^{T ( h ) } \| A (h)  e^{ - i t P (h) / h } \chi( P (h) 
)u\|^2 _{\HH_1 } dt \,.
\end{equation}
\end{thm}

To motivate the abstract presentation we relate the notation
of Theorem \ref{t:3} to a concrete situation.
Thus let $ P ( h ) = -h^2 \Delta $ be the Dirichlet Laplacian 
on a compact manifold $ \Omega $, with boundary $ \partial \Omega $.
Then 
\[ \HH = L^2 (\Omega ) \,, \ \ {\mathcal D} = H^{2}(\Omega )\cap 
H^1_0 ( \Omega ) \,.\]
Let $ \Gamma \subset \Omega $. We then define 
\[ \HH_1 = L^2 (\Gamma ) \,, \ \ {\mathcal D} 
\ni u \longmapsto A ( h )u  = h \partial_\nu u \rest_\Gamma \in \HH_1 \,, \]
where $ \partial_\nu $ denotes the inward pointing normal to $ \partial
\Omega $. The estimate \eqref{eq:3.3} is a typical {\em observability 
estimate} equivalent by duality to an {\em exact control} 
statement (see Sect.\ref{geomcont}). An abstract method for obtaining semi-classical estimates 
\eqref{eq:3.1} will be presented in Sect.\ref{sbb}.

\begin{proof}
Let us put $ v (t) = \exp( - it P ( h )/h ) \chi ( P ( h )) 
u $. We introduce a function
$ \psi \in \CIc ( \RR ; [0,1] )$, 
and put
$$ w ( t ) = \psi \left( \frac{ t }{ T ( h )} \right) v ( t) \,.$$
Clearly,
\[ ( i h \partial_t -  P ) w ( t) = \frac{ i h }{ T( h ) } \psi' \left( 
\frac{t}{ T( h )} \right) v ( t ) \,.\]
Because of the compact support we can take the (semi-classical) Fourier
transform in $ t $ which gives
\[ ( \tau + P ) \widehat w ( \tau ) = - \frac{ i h }{ T( h ) } 
{\mathcal F}_{{t\rightarrow\tau}}( \psi' ( \bullet / T ( h ) ) v ) ( \tau ) \,.\]
For $ \tau \in I $ we can use \eqref{eq:3.1} which gives
\[ \| \widehat w ( \tau ) \|_{\HH} \leq 
\frac{ G ( h )}{ T ( h ) } \| {\mathcal F}_{t\rightarrow\tau}( 
 \psi' ( \bullet / T ( h ) ) v )  
( \tau )  \|_\HH + g ( h ) \| A (h) \widehat w ( \tau ) \|_{\HH_1} \,.\]
Using the generalized Plancherel theorem we obtain
\[ \int_{I} \| \widehat w ( \tau ) \|^2_{\HH} d \tau \leq 
2 \frac{ G ( h )^2}{ T ( h )^2 } \| \psi' ( \bullet/ T( h ) )  v   
  \|_{ L^2 ( \RR_t ;\HH)} ^2  + 2 
g ( h )^2  \| A (h)  w  \|^2_{ L^2 ( \RR_t ; \HH_1 )} \,.\]
We now want to show that we can integrate over
$ \RR $ in place of $ I $ in the left hand side. That follows from 
\begin{equation}
\label{eq:3.lem} 
\| \hat w ( \tau ) \bbbone_{ \RR \setminus I }( \tau )  \|_\HH 
= {\mathcal O}\left( \left(\frac{h}{1 + |\tau | } \right)^\infty \right)
\| \chi ( P ) u \|_\HH \,,
\end{equation}
which in turn follows from integration by parts in 
\[ \hat w ( \tau ) = \int_\RR e^{ -it ( P + \tau)/h } \psi \left(\frac{t}{T} \right) 
\chi ( P ) u dt = 
\int_\RR ( - ( P + \tau ) ^{-1} hD_t  e^{ -it ( P + \tau)/h } ) \psi \left(\frac{t}{T} \right) 
\chi ( P ) u dt \,, \]
using 
\[  \forall \; \tau \in \RR \setminus I \ \  \| ( P + \tau ) \chi ( P ) u 
\|_\HH  \geq \frac1 C \| \chi ( P ) u \|_\HH \,.\]

Thus we obtained
\[ \| w \|^2_{ L^2 ( \RR_t ; \HH )} \leq 2\frac{ G( h )^2 } { T ( h )^2 }
\| \psi' ( \bullet/ T ( h ) ) v \|^2_{ L^2 ( \RR_t ; \HH ) } 
+ g ( h )^2  \| A (h)  w  \|^2_{ L^2 ( \RR_t ; \HH_1 )} 
+ {\mathcal O} ( h ^\infty ) \| \chi (P ) u \|^2 \,,\]
and the first term on the right can be absorbed on the left 
using \eqref{eq:3.2}. In fact, since
\[ \sup_{ \phi \in \CIc ( ( 0 , 1 )) } \frac{ \int_0^1 \phi (s)^2 ds
}{ \int_0^1 \phi' (s)^2 ds } = \pi^{-2}\,, \]
we have from the definition of $ w $, and for any $ \epsilon > 0 $,
\[  \| \chi( P ) u \|_{\HH } \leq 2  (\pi^2 + \epsilon)
\frac{ G( h )^2 } { T ( h )^2 } \| \chi ( P )  u \|_{\HH} 
+ 2 \frac{g ( h )^2}{ T(h) }   \| A (h)  w  \|^2_{ L^2 ( \RR_t ; \HH_1 )} 
+ {\mathcal O} ( h ^\infty ) \| \chi (P ) u \|_{\HH}^2 \,. \]
This completes the proof once we take $ h $ small enough.
\end{proof}

\section{Semiclassical black box resolvent estimates}
\label{sbb}
In this section we will make assumptions under which 
resolvent estimates can be obtained in the semi-classical 
setting. For simplicity no boundary will be allowed here.

Let $ X $ be a compact $
\CI $  manifold. Let
$ P(h) \in \Psi_h^{2,0} ( X ; \Omega_X^{\frac12} ) $ 
be formally self-adjoint on $ L^2 ( X ; \Omega_X^{\frac12} ) $. 
We assume that, if $ p $ is the principal symbol of $ P ( h ) $
then 
\begin{equation}
\label{eq:4.1}  p = 0 \Longrightarrow dp \neq 0 \,,  
\ \   p \geq \langle  \xi \rangle^2 /C  \ \text{ for $ |\xi| \geq C $,} \end{equation}
and that for some $ \delta > 0 $
\begin{equation}
\label{eq:4.2}
 p^{-1} ( [ - \delta, \delta ] ) \Subset T^* X \,.
\end{equation}
\begin{figure}[ht]
$$\ecriture{\includegraphics[width=12cm]{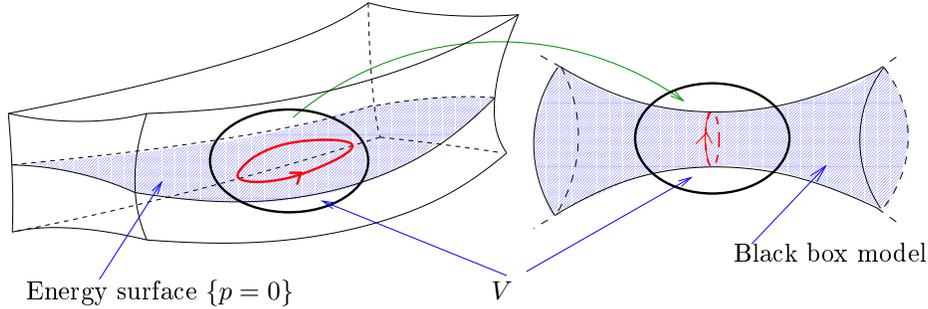}}
{\aat{2}{-1}{Energy surface $\{p=0\}$}\aat{27}{-1}{$V$}
\aat{40}{1}{Black box model}}$$
\caption{A semi-classical black box with an hyperbolic trapped trajectory.}
\label{fig:sbb}
\end{figure}

Suppose that $ Q (h ) $ is a family of 
bounded operators on a Hilbert 
space $ \HH $. Suppose that there exist bounded operators
\[\begin{split}
&  U_1 (h) \; : \;   L^2 ( X ; \Omega_X^{\frac12} ) \ \longrightarrow \HH \\
& U_2 (h ) \; : \; \HH \longrightarrow L^2 ( X ; \Omega_X^{\frac12} ) \,, \\
& \chi^\sharp (h) \; : \; \HH \longrightarrow \HH \,, 
\end{split}
\]
such that, microlocally near $ V $, an open subset of $ p^{-1} ([-\delta 
, \delta ] ) $, we have  
\begin{equation}
\label{eq:4.0}
\begin{split}
& U_2 ( h ) \circ U_1 ( h ) = Id \,, \\  
& U_1 ( h ) \circ U_2 ( h  ) = \chi^\sharp ( h ) \,, \\ 
& U_1 ( h ) \circ P ( h ) \circ U_2 ( h ) 
 = Q( h ) \circ \chi^\sharp ( h )  \,. \end{split} \end{equation}
In practice, the operators $ U_j ( h ) $ are $h$-Fourier integral operators
(see Proposition \ref{p:2.0}) but we do not need to make this assumption
in the abstract presentation. Figure \ref{fig:sbb} shows our setup schematically 
in the case relevant for the proof of Theorem \ref{t:cdvp}.
\begin{thm}
\label{t:4}
Let $ P( h ) $ and $ Q ( h ) $ satisfy the assumptions above and $V_{0}$ be an open relatively compact subset of $T^*X$.
Suppose that $ A \in \Psi^{0,0}_h ( X , \Omega_X^{\frac12} ) $ 
is microlocally elliptic in $ V_0 $, 
and that there exists $ T > 0 $ such that 
\begin{gather}
\label{eq:4.3}
\begin{gathered} 
\forall \; \rho \in  p^{-1} ( 0 ) \setminus V \; \exists \ 0 < t < T \,, \ \epsilon \in 
\{\pm 1\} \\ \exp ( \epsilon s H_p ) ( \rho ) \subset p^{-1} ( 0 ) \setminus V \,, \ 0<s<t \,, \ \ 
\exp ( \epsilon t H_p ) ( \rho ) \in V_0 \,.
\end{gathered}
\end{gather}
Suppose also that
\begin{equation}
\label{eq:4.4}
\| \chi^\sharp ( h )  Q( h ) ^{-1} \| \leq \frac{ G ( h ) }{h }, \, G(h)\geq 1 \,.\end{equation}
Then for $ u \in \CI ( X , \Omega_X^{\frac12})  $ we have 
\begin{equation}
\label{eq:4.5}  \| u \| \leq C \frac{ G ( h ) }h \| f \| + 
G ( h ) \| A u \| \,. \end{equation}
\end{thm}
We start with the following standard: 
\begin{lem}
\label{l:4}
Suppose that $ p,A $, and $ V $ satisfy \eqref{eq:4.3}. 
If  $ B \in \Psi^{0,0} ( X , \Omega_X^{\frac12} ) $ and
$ WF ( B ) \subset T^* X \setminus V $ then
\begin{equation}
\label{eq:4.6} 
\| B u \| \leq C h^{-1}\| P u \| + \| A u \| + {\mathcal O}(h^\infty) \|u\|
\,.
\end{equation}
\end{lem}
\begin{proof}
In view of  the compactness of
$ p^{-1}(0) $ we can replace $ V_0 $ by a precompact neighbourhood
of $ V_0 \cap p^{-1} ( 0 ) $.
The assumption  \eqref{eq:4.3} then shows that
it is enough to prove a local version of the estimate. We 
can suppose that $ WF ( A ) \subset U $ where $ U$ is a small
neighbourhood of $ m_0 \in V_0 $ and 
$$ WF ( B ) \subset 
\bigcup_{ 0 \leq t \leq t_0 } \exp( \epsilon t H_p ) ( U_1 ) \subset 
T^* X \setminus V \,, \ \ U_1 \Subset U \,. $$
If $ t_0 $ is small enough we can apply Proposition \ref{p:2.1},
as the estimate is clear in the case of $ P = hD_{x_1} $. 
In general, we can then split the interval $ [0, t_0 ] $ 
into subintervals in which the $t_0$-small argument can be applied.
\end{proof}

\noindent
{\em Proof of Theorem \ref{t:4}.}
Suppose that $ B_1 $ satisfies
\[ WF ( B_1 ) \subset V_1 \,, \ \ V \Subset V_1 \,,  \  \
WF ( I - B_1 ) \subset T^* X \setminus V \,. \]
Then if $ V_1 $ is sufficiently close to $ V $ then
using the second part of \eqref{eq:4.0} we have 
\begin{equation}
\label{eq:4.nnew} \begin{split}
\| B_1 u \| & = \| U_2 \chi^\sharp U_1 B_1 u \| + {\mathcal O} ( h ^\infty ) \|u\|
\\
& = \| U_2 \chi^\sharp Q^{-1} Q U_1 B_1 u \|  + {\mathcal O} ( h ^\infty )  \| u \|
\end{split} \end{equation}
If we now apply \eqref{eq:4.4} and then \eqref{eq:4.0} again, we 
obtain
\begin{equation}
\label{eq:4.new}
 \begin{split}
\| B_1 u \| & \leq   \frac{G( h )}{ h}  \| Q \chi^\sharp
U_1 B_1u \|_\HH + {\mathcal O} ( h ^\infty )  \|u\| \\
& \leq   C \frac{G( h )}{ h}  (  \| P u \|  + \| [ P, B_1 ] u \| )
+ {\mathcal O} ( h ^\infty )  \|u\| \\
& \leq    C \frac{G( h )}{ h}    \| P u \|  +  G ( h ) \| B_2 u \|  
+ {\mathcal O} ( h ^\infty )  \|u\| \,, 
\end{split} \end{equation}
where  $ B_2 \in \Psi^{0,0} ( X , \Omega_X^{\frac12} ) $ satisfies
\[  WF ( B_2 ) \subset V_1 \setminus V \,, \ \ WF ( (I - B_2 ) [ P, B_1] )
= \emptyset \,.\]
Lemma \ref{l:4} now shows that 
\[ \| B_1 u \| \leq  C \frac{G( h )}{ h} \| P u \| + G ( h ) \| A u \| +
{\mathcal O} ( h ^\infty ) \| u \| \,.\]

We now choose $ B_3 \in \Psi_h^{0,0} ( X , 
\Omega_X^{\frac12} ) $ such that $ WF ( B_3 ) \subset T^* X \setminus
V $, and $ WF( I - B_3 ) \subset V_1 $. We can apply Lemma \ref{l:4} 
with $ B = B_3 $ and that gives \eqref{eq:4.5} as $ \| u \| 
\simeq \| B_1 u \| + \| B_3 u \| $.
\stopthm

In some situations we can obtain improved estimates under
a modified assumption on $ Q^{-1}$. This modification will be
crucial in Sect.\ref{eaa} where we will prove \eqref{eq:1.cdv}.
We present it separately not to obscure the simplicity of 
Theorem \ref{t:4}:

\medskip
\noindent
{\bf Theorem \ref{t:4}$'$.}\hspace{0.1in}{\em Suppose that the assumptions of
Theorem \ref{t:4} hold, and that in addition,
\begin{equation}
\label{eq:4.7} \| \chi^\sharp ( h ) Q^{-1} U_1 \phi ( h ) \| \leq \frac{g(h)}{h} 
\,, 
\end{equation}
where $ \phi ( h ) $ is a microlocal cut-off to a neighbourhood of
$ V_1 \setminus V $, where $ V_1 \Supset V $ is a small neighbourhood
of $ V$. 
Then we have,
\begin{equation}
\label{eq:4.8}
 \| u \| \leq C \frac{ G ( h ) }h \| Pu  \| + 
g ( h ) \| A u \| \,. \end{equation}
}

\medskip

\begin{proof} We revisit the proof of Theorem \ref{t:4}. Instead of
moving instantly to \eqref{eq:4.new} from \eqref{eq:4.nnew} using
\eqref{eq:4.5}, we apply the identities~\eqref{eq:4.0}, and write
\[  \begin{split}
\| B_1 u \| & \leq  C   \| \chi^\sharp Q^{-1} Q 
\chi^\sharp  U_1 B_1u \|_\HH + {\mathcal O} ( h ^\infty )  \|u\| \\
& =  \| \chi^\sharp Q^{-1} U_1 ( B_1 P u + [ P, B_1] u ) \|_\HH
+ {\mathcal O} ( h ^\infty )  \|u\| \\
& \leq   \|\chi^\sharp Q^{-1} U_1  B_1 P u \|_\HH + 
   \|\chi^\sharp Q^{-1} U_1 \phi ( h ) [ P, B_1] u  \|_\HH 
+ {\mathcal O} ( h ^\infty )  \|u\| \,, 
\end{split} \]
where we could insert the cut-off $ \phi(h) $ due to the microsupport
properties of $ B_1$. 

If we apply \eqref{eq:4.5} and \eqref{eq:4.7} we obtain a local 
version of \eqref{eq:4.8}:
\[ \| B_1 u \|  \leq  C   \frac{ G ( h ) }h \| Pu  \| + 
\frac{ g (h) }{ h} \| [ P, B_1 ] u \| \,.\]
The proof is then completed as in the case of Theorem \ref{t:4}.
\end{proof}

\section{Estimates in the homogeneous case: classical control}
\label{cbb}

In this section we will adapt the semi-classical arguments of Sect.\ref{sbb}
to obtain a classical version of the estimate \eqref{eq:4.5}. 
We start by modifying the black box assumptions where we essentially
follow \cite{SjZw91},\cite{Sj97} but change the ambient space from 
$ \RR^n $ to an arbitrary manifold.

Thus let $ X $ be compact $ \CI $ manifold with a (possibly empty) 
boundary $ \partial X $. We consider an elliptic differential
operator of order two, 
$$ P_0 \in {\rm Diff}^2 ( X , \Omega_X^{\frac12} ) \,, $$
with a domain $ {\mathcal D}_0 \subset L^2 ( X , \Omega_X^{\frac12} )$.
The choice of the domain includes the possible boundary conditions.

\begin{figure}[ht]
$$\ecriture{\includegraphics[width=12cm]{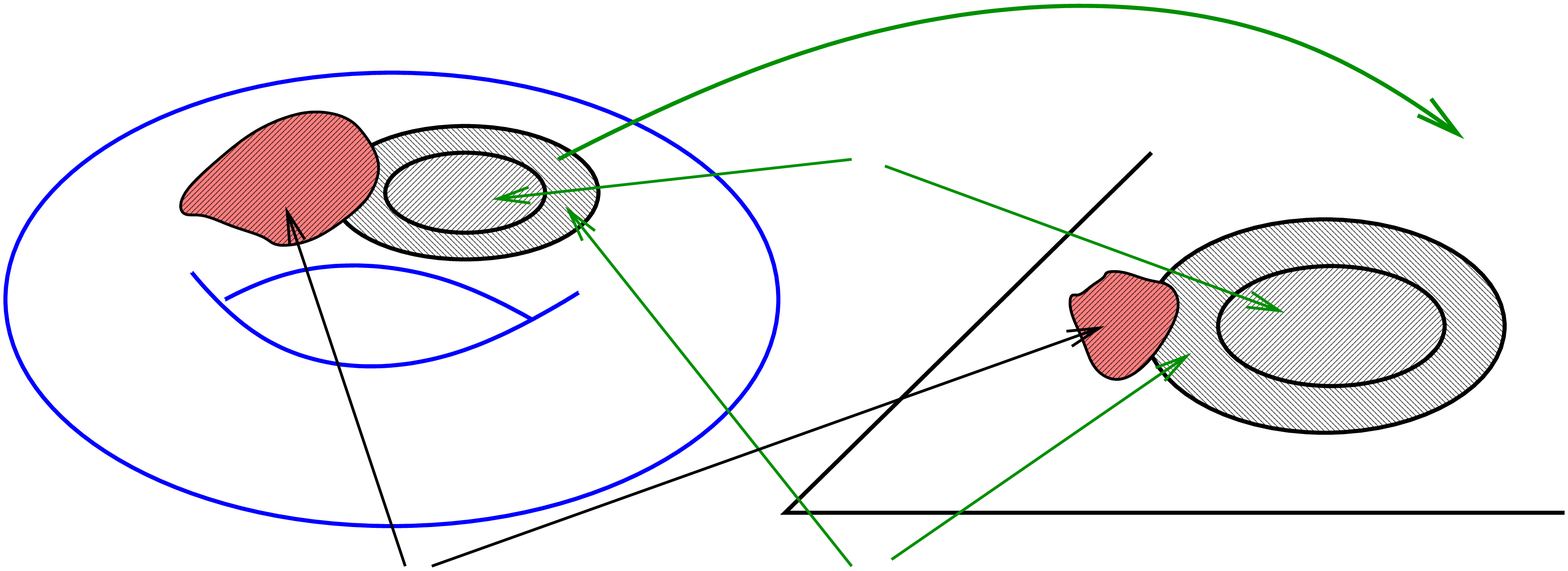}}
{\aat{10}{-1}{Control regions}
\aat{38}{-1}{Black box model}
\aat{27}{-1}{$\widetilde Y$}\aat{27}{14}{$Y$}\aat{27}{2}{$\mathbb{R}^2$}
\aat{2}{7}{$\mathbb{T}^2$}
}$$
\caption{The black box $ Y$, its neighbourhood, $ \widetilde Y $, 
in the case when $ X = {\mathbb T}^2 
$ is the flat torus, and $ \widetilde X = \RR^2 $, the plane.}
\label{fig:bb}
\end{figure}

Let $ Y \subset X $ be an open set. We also consider an auxiliary manifold $ \widetilde X $, which
coincides with $ X $ on a neighbourhood, $ \widetilde Y $ of $ Y $ -- 
see Fig.\ref{fig:bb} for a visualization.

We then consider complex Hilbert spaces $\mathcal{H}$, $ \HH_{\rm{bb}} $
with orthogonal
decompositions
\[ \begin{split}
& {\mathcal{H}}={\mathcal{H}}_{Y} \oplus  L^2(
X \setminus Y , \Omega_X^{\frac12} )\\
& 
{\mathcal{H}}_{\rm{bb}} ={\mathcal{H}}_{Y} \oplus 
L^2(
\widetilde X \setminus Y , \Omega_{\widetilde X}^{\frac12} )\,.
\end{split} \]
For $ \HH$  the orthogonal projections on the two factors are denoted by 
$ \bbbone_{ Y } $  and $ \bbbone_{ X\setminus Y} $
respectively. If $ \chi_j \in \CI ( X ) $ satisfy
\begin{equation}
\label{eq:chi}
\supp \chi_0 \subset \complement \supp ( 1 - \chi_1 ) \subset 
\supp \chi_1 \subset \widetilde  Y \,, \ \ 
\supp ( 1 - \chi_0 ) \subset X \setminus \widetilde Y 
\end{equation}
then multiplication by $ \chi_j $ is well defined on $ \HH $ and $ \HH_{\rm{bb}} $.

On $ L^2 ( X )  $ and $ \HH_{\rm{bb}} $ we have unbounded operators,
$ P_0 $ and $ P_{\rm{bb}} $ respectively with domains
\[ \begin{split}
& {\mathcal D}_0 \stackrel{\rm{def}}{=} {\mathcal D} ( P_0) \subset 
L^2 ( X , \Omega_X^{\frac12} ) \\
& {\mathcal D}_{\rm{bb}} \stackrel{\rm{def}}{=} {\mathcal D} ( P_{\rm{bb}})
 \subset  \HH_{\rm{bb}} \,. \end{split} \]

A self-adjoint operator, $ P :\mathcal{H} \longrightarrow \mathcal{H}$,
has the  domain $\mathcal{D} \subset \mathcal{H}$, satisfying the following
conditions:
\begin{gather*}
 \bbbone_{X \setminus Y } {\mathcal{D}} =
 \bbbone_{X \setminus Y }  {\mathcal D_0 }  \,, \ \ 
 \bbbone_{Y } {\mathcal{D}} =
 \bbbone_{Y }  {\mathcal D} _{\rm{bb}}  \,, \ \ 
\\
( 1 - \chi_1) P = ( 1 - \chi_1 )P(1-\chi_{0})= (1-\chi_{1}) P_0 ( 1 - \chi_0 )= (1-\chi_{1}) P_0  \,, \\
P \chi_{0}= \chi_1 P\chi_{0} = \chi_1 P_{\rm{bb}} \chi_0= P_{\rm{bb}} \chi_0  \,,
\end{gather*}
for any functions satisfying \eqref{eq:chi}.
We use the notation from \cite{SjZw91} and in particular write
\[  {\mathcal D}^\infty = \bigcap_{ k \in \NN } {\mathcal D } ( P^k ) \,, \ \ 
 {\mathcal {D}}^{-N}= ({\mathcal {D}}^N)^* \,.\]
We also make another standard ``black box'' assumption:
\[ (P+i)^{-1}  \ \text{ is compact on ${\mathcal {H}}$.}\]

As in previous sections we have two types of results. To obtain 
the assumptions of an analogue of Theorem \ref{t:3} we need resolvent 
estimates based on {\em black box resolvent estimates}. That is
provided in 
\begin{thm}
\label{t:5.1}
Suppose that $  A : {\mathcal D} (A ) \rightarrow \HH_1 $, 
$ A_{\rm{bb}}
 : {\mathcal D} (A_{\rm{bb}}) \rightarrow \HH_1 $, where $ \HH_1 $ is
a Hilbert space, 
$  {\mathcal D} (A) \supset {\mathcal D}^\infty $, 
$ {\mathcal D} (A_{\rm{bb}}) \supset {\mathcal D}^\infty_{\rm{bb}} $,
satisfy, for $ u \in {\mathcal D}^\infty $ and $ v \in {\mathcal D}_{\rm{bb}}
^\infty $,
\begin{gather}
\label{eq:blackboxcontrol}
\begin{gathered}
 \| \bbbone_{X \setminus Y } u \|_\HH \leq C \langle \lambda \rangle^{-\frac12}
\| ( P - \lambda ) u \|_{\mathcal H} + \| A u \|_{\HH_1} + {\mathcal O} (
\langle \lambda \rangle^{-\infty}) \| u \|_{\HH} \,,\\
 \| v \|_{\HH_{\rm{bb}}} \leq G( \lambda ) \left( 
\langle  \lambda \rangle^{-\frac12}  \| ( P_{\rm{bb}}
 - \lambda )v \|_{\HH_{\rm{bb}}}
+ \| A_{\rm{bb}}  v \|_{\HH_1}\right)  \,, \ \
|\lambda | \rightarrow \infty \,, \\ 
   A\chi_0  =   \chi_1 A \chi_0 =  \chi_1 A_{\rm{bb}} \chi_0= A_{\rm{bb}}\chi_{0}  \,, \ \ 
G ( \lambda ) \geq 1 \,,\\
\forall k \exists C_{k}; \|A_{bb}\chi_{0}u\| \leq C_{k}( \|Au\|_{\HH_{1}}+ \|u\|_{{\mathcal {D}}^{-k}})
\end{gathered}
\end{gather}
for any $ \chi_j$'s satisfying \eqref{eq:chi}.
Then
\begin{equation}
\label{eq:5.0}
\| u \|_{\HH} \leq C_1 G ( \lambda) \left( \langle \lambda \rangle^{-1/2}\| ( P - \lambda ) u  \|_{\HH} 
+  \| A u \|_{\HH_1} \right) \,.
\end{equation}
\end{thm}
\begin{proof}
We first prove the following estimate:
\begin{equation}\label{eq4.500} \langle \lambda\rangle ^{-1/2}\| [ P_{0}, \chi_{0}]u
\|_{\HH}\leq C\left(\langle \lambda\rangle ^{-1/2}
\| ( P - \lambda ) \|_{\HH} + \|u\bbbone_{X \setminus Y }\|_{\HH}\right)
\end{equation}
Indeed, the ellipticity of $P_{0}$ gives
\begin{equation}
\begin{split}
\|u\|_{H^2( \text{supp}( \nabla \chi_{0}))}
& \leq C \left( \|P_{0}u\|_{L^2( X\setminus Y)}+\|u\|_{L^2( X\setminus Y)}\right)\\
& \leq  C \left( \|(P-\lambda)u\|_{L^2( X\setminus Y)}+
\langle \lambda\rangle\|u\|_{L^2( X\setminus Y)}\right) 
\end{split}
\end{equation}
Using the inequality
$\|u\|_{H^1}\leq C \sqrt{ \|u\|_{H^2}\|u\|_{L^2}}$ we get~\eqref{eq4.500}.

 We now turn to the proof of~\eqref{eq:5.0}. The black box assumptions give 
$$(P_{\rm bb}-\lambda) \chi_{0}u=(P-\lambda) \chi_{0}u= [P_{0}, \chi_{0}]u + \chi_{0}(P-\lambda)u $$
Using~\eqref{eq:blackboxcontrol} we obtain
$$\|\chi_{0}u\|_{\HH}=\|\chi_{0}u\|_{\HH_{\rm bb}}\leq G(\lambda) \left(\langle \lambda \rangle^{-1/2}(\| \chi_{0}(P-\lambda)u\|_{\HH_{\rm bb}}+\|[P_{0}, \chi_{0}]u\|_{\HH_{\rm bb}}) +  \| A_{\rm bb}\chi_{0}u \|_{\HH_{1}} \right)$$
Above, we can replace the norms in $\HH_{\rm bb}$ by norms in $\HH$ and,
using ~\eqref{eq:blackboxcontrol} and~\eqref{eq4.500}, this implies
$$\|\chi_{0}u\|_{\HH}\leq CG(\lambda) \left(\langle \lambda \rangle^{-1/2}\|(P-\lambda)u\|_{\HH}+\|Au\|_{\HH_{1}}+ {\mathcal {O}}( \langle \lambda\rangle^{-\infty})\|u \|_{\HH} + \|u\|_{{\mathcal {D}}^{-k}}\right)$$ 
To conclude the proof 
we use the first inequality in ~\eqref{eq:blackboxcontrol} and the fact that
$$\|u\|_{{\mathcal {D}}^{-k}} = \|(P+i)^{-k}u\|_{\HH}\leq C_{k}\left(\langle \lambda\rangle ^{-k}\|u\|_{\HH}+ \| ( P - \lambda ) u \|_{\HH}\right)$$
\end{proof}
\begin{rem}
In the proof above, the operators $A$ and $A_{\rm bb}$ could depend on 
$\lambda$ as long as the assumptions are uniform in $\lambda$.
\end{rem}

The difference between the semi-classical and classical control estimates,
\eqref{eq:3.3} and \eqref{eq:5.2} below, is more serious. In the classical
case the low energy contribution does not allow an explicit time dependent
constant we have in \eqref{eq:3.3} (compare \eqref{eq:5.2} and \eqref{eq:5.new}
below). As investigated recently in \cite{Mi03} violent behaviour is
expected when fast control is a goal.

\begin{thm}
\label{t:5.2}
Suppose that $ A : {\mathcal D} ( A)  \rightarrow \HH_1 $, where $ \HH_1 $ is
a Hilbert space, 
$  {\mathcal D} (A) \supset {\mathcal D}^\infty $, satisfies the following condition:
for all $ N $ there exists $ C_N $ such that for all $ k \in \NN $ and 
$ u \in {\mathcal D}^\infty $,
\begin{equation}
\label{eq:cut}
   \| A \psi (2^{-k} P ) u \|_{\HH_1} + \| A ( 1- \psi )(2^{-k}P) u \|_{\HH_1}
 \leq 
C_{N}\left(\| A u \|_{\HH_1}+ 2^{-kN}\|u\|_{{\mathcal {D}}^{-N}} \right)\,,
\end{equation}
$ \psi \in C^\infty_{0}( (1/2, 2))$.

Suppose also that
for all $ \lambda \in \RR $ and $ u \in {\mathcal D}^\infty $
we have
\begin{equation}
\label{eq:5.1}
 \| u \|_{\HH} \leq G ( \lambda ) \| ( P - \lambda ) u \|_\HH 
+ g ( \lambda ) \| A u \|_{\HH_1} \,, \end{equation}
where $ G $ and  $ g $ satisfy
\begin{equation}\label{eq:5.66}
\langle \lambda \rangle^{-1} \leq G ( \lambda ) \leq \langle \lambda \rangle^{N_0} \,, \ \ 
C\langle \lambda \rangle ^{-N_{0}}\leq g(\lambda)\leq C'\langle \lambda \rangle ^{N_{0}}\,, 
\ \ g(\lambda/2)\leq C g(\lambda) \leq C' g(2\lambda).
\end{equation}
We also assume the following weak continuity property (see Remark~\ref{re:5.5} for a discussion)
There exist $ N_1 \in \NN $ and a Hilbert space ${\mathcal {H}}^ \sharp$ such that
$ {\mathcal {H}}_{1}\subset {\mathcal {H}}^\sharp $ continuously, 
and the operator $Ae^{itP}$ is continuous from
${\mathcal {D}}^{-N_{0}}$ to 
$H^{-N_1}_{\text{loc}}( {\mathbb R}_{t}; {\mathcal {H}}^\sharp)$.
\par

Then there exist constants $ C_0 $ and $ C_1=C_1 (T)  $ 
such that for {\em any } $ T > C_1 \limsup_{|\lambda| \rightarrow \infty }
G ( \lambda ) $ we have 
for $ u \in {\mathcal D}^\infty   $, 
\begin{equation}
\label{eq:5.2}
\| \langle g ( P ) \rangle^{-1} u \|^2_\HH
 \leq C_1 (T)  \int_0^T \| e^{it P} A u \|_{\HH_1} ^2 dt \,.
\end{equation}
\end{thm}
\medskip
\noindent
\begin{rem}\label{re:5.5}
In the case where the operator $P$ is the Laplace operator with Dirichlet boundary conditions, the weak continuity is satisfied in the two following typical  situations:
\begin{enumerate}
\item If $A$ is a pseudodifferential operator supported in the interior of $X$, then ${\mathcal {H}}_{1}= L^2(X)$ and we can take ${\mathcal {H}}^\sharp$ to be another Sobolev space $H^{-s}(X)$.  
\item If $Au= \partial_{n}u \rest_{\Gamma}$ where $\Gamma\subset \partial \Omega$ and $\partial _{n}$ is the normal derivative to the boundary. Then we can take 
${\mathcal {H}}^\sharp  = {\mathcal {H}}_{1}= L^2(\partial X)$ as standard trace regularity results 
for solutions of Schr\"odinger equations show that the assumptions hold with $ N_1 $ sufficiently 
large.
\end{enumerate}
\end{rem}

\noindent
{\em Proof of Theorem~\ref{t:5.2}:}
We follow closely the proof of Theorem \ref{t:3} observing first that, with $\Psi\in C^\infty_{0}(]1/2, 2[)$ equal to $1$ close to $1$,  
\eqref{eq:5.1} and \eqref{eq:cut} imply
\begin{equation}
\langle g ( \lambda ) \rangle^{-1} \|\Psi ( P / 
{\langle \lambda \rangle}) u \|_\HH \leq
G ( \lambda )\langle g ( \lambda ) \rangle^{-1} \| \Psi ( P /{\langle \lambda \rangle}) ( P - \lambda ) u \|_\HH 
+ C \| A u \|_{\HH_1}+ C \langle \lambda \rangle ^{-N}\|u\|_{{\mathcal {D}}^{-N}}\end{equation}
which in turn implies 
\begin{equation}
\|\langle g ( P ) \rangle^{-1}\Psi ( P/
 {\langle \lambda \rangle} ) u \|_\HH \leq
G ( \lambda ) \| \langle g ( P ) \rangle^{-1} ( P - \lambda ) u \|_\HH 
+ C\| A u \|_{\HH_1}+ C \langle \lambda \rangle ^{-N}\|u\|_{{\mathcal {D}}^{-N}} \,.
\end{equation}
The functional calculus of self adjoint operators gives 
\[ \|(1-\Psi) ( P /{\langle \lambda \rangle})\langle g ( P ) \rangle^{-1} u\|_\HH \leq \sup_{\xi}\left| \frac {(1- \Psi)( \xi /{\langle \lambda \rangle})}{\xi - \lambda}\right|\|\langle g(P) \rangle ^{-1}u \|_{\HH}\leq \frac C {1+ |\lambda|}\|\langle g(P) \rangle ^{-1}u \|_{\HH}\,,\]
which, using \eqref{eq:cut} again, and~\eqref{eq:5.66} implies (taking $N$ large enough) that for $|\lambda|$ large enough, 
\begin{equation}\label{eq:5.60} \|\langle g ( P ) \rangle^{-1} u \|_\HH \leq
CG ( \lambda ) \| \langle g ( P ) \rangle^{-1} ( P - \lambda ) u \|_\HH 
+ C\| A u \|_{\HH_1}\, .
\end{equation}
Proceeding as in the proof of Theorem~\ref{t:3} we define $ v (t) = \exp(  it P  )
u $. We introduce a function
$ \psi \in C^\infty_{0} ( ]0,1[ )$, 
and put
$$ w ( t ) = \psi \left( \frac{ t }{ T } \right) v ( t) \,,$$
so that
\[ ( i  \partial_t -  P ) w ( t) = \frac{ i  }{ T } \psi' \left( 
\frac{t}{ T} \right) v ( t ) \,.\]
Because of the compact support we can take the  Fourier
transform in $ t $ which gives
\[ (  \tau  - P ) \widehat w ( \tau ) = \frac{ i  }{ T} 
{\mathcal F}_{{t\rightarrow\tau}}( \psi' ( \bullet / T  ) v ) ( \tau ) \,.\]  

Let $\varrho$ be a large constant to be fixed later. For $\langle \tau \rangle \geq \varrho/2$ we estimate $ \widehat w ( \tau )$ using~\eqref{eq:5.60} which gives
\begin{equation}\label{eq:5.63}
\| \langle g( P)\rangle^{-1}\widehat w ( \tau )\|_{{\mathcal {H}}}\leq
C\frac{G ( \tau )} T \| \langle g ( P ) \rangle^{-1} {\mathcal F}_{{t\rightarrow\tau}}( \psi' ( \bullet / T  ) v ) ( \tau )\|_\HH 
+ C\| A \widehat w ( \tau ) \|_{\HH_1}\, ,
\end{equation}
For $\langle \tau \rangle \leq \varrho/2$ we simply write, with $\chi\in C^\infty_{0}( ]-1, 1[)$ equal to $1$ on $[-1/2, 1/2]$,
\begin{equation}\label{eq:5.62}
\widehat w ( \tau )= \int_{t\in {\mathbb R}}e^ {it(P- \tau)} \psi \left( \frac{ t }{ T } \right) (\chi(P/\varrho) u+ (1- \chi(P/\varrho)) u) dt
\end{equation}
The contribution of the first term is bounded (in ${\mathcal {H}}$) by $\| (\chi(P/\varrho) u\|_{\HH}$ and by integrations by parts with the operator $\frac {i\partial _{t}}{ P-\tau}$ we can bound the contribution of the second term by 
\begin{equation}\label{eq:5.61}
C_{N}\| \frac 1 { (1+ |T|+ |\varrho| + \langle P \rangle)^N }u\|_{{\mathcal {H}}}
\end{equation}
From~\eqref{eq:5.63},\eqref{eq:5.62},\eqref{eq:5.61} and the bounds on the weight $g$, we get
\begin{multline}
\| \langle g( P)\rangle^{-1}\widehat w ( \tau )\|^2_{L^2( {\mathbb R}_{\tau}; {\mathcal {H}})}\leq
C \left( \frac{ \sup_{|\tau| \geq \rho/2 } G ( \tau)} {T} \right)^2 
\| \langle g ( P ) \rangle^{-1} {\mathcal F}_{{t\rightarrow\tau}}( \psi' ( \bullet / T  ) v\|^2 _{L^2( {\mathbb R}_{\tau}; \HH)}\\
 +  C\int_0^T \| Ae^{  i t P } 
 u \|_{\HH_1}^2 dt + C\| \bbbone_{ \langle
 P \rangle \leq \rho } u \|_{\HH}^2 + C \langle \varrho \rangle ^{-N_{0}} \|u \|^2_{{\mathcal {D}}^{-N_{0}}}\end{multline}
Remark that 
\begin{equation}
\| \langle g( P)\rangle^{-1}\widehat w ( \tau )\|_{L^2( {\mathbb R}_{\tau}; {\mathcal {H}})}= T^{1/2 }\|\Psi\|_{L^2}\|\langle g( P)\rangle^{-1}u\|_{\HH} 
\end{equation}
and
\begin{equation}
\| \langle g( P)\rangle^{-1}{\mathcal F}_{{t\rightarrow\tau}}( \psi' ( \bullet / T  ) v\|_{L^2( {\mathbb R}_{\tau}; {\mathcal {H}})}= T^{1/2 }\|\Psi'\|_{L^2}\|\langle g( P)\rangle^{-1}u\|_{\HH} 
\end{equation}
Consequently, taking $\rho$ large enough the assumption $ T > C_1 \limsup_{|\lambda| 
\rightarrow \infty } G ( \lambda ) $ ensures that we can eliminate the first and the last  terms in the right hand side  and get
\begin{equation}\label{eq5.new}
\| \langle g( P)\rangle^{-1}u\|_{\HH}\leq C \int_0^T \| Ae^{  i t P } 
 u \|_{\HH_1}^2 dt + C\| \bbbone_{ \langle
 P \rangle \leq \rho } u \|_{\HH}^2
\end{equation}
 To eliminate the last term we use the compactness-uniqueness argument from~\cite{BLR} which we now
recall. Proceeding by contradiction we obtain a sequence $(u_{n})$ such that
\begin{equation}
C \| \bbbone_{ \langle
 P \rangle \leq \rho } u_{n} \|_{\HH}^2\geq 1= \| \langle g ( P ) \rangle^{-1} u_{n} \|^2 _\HH \geq n \int_0^T \|A e^{  i t P } 
 u_{n} \|_{\HH_1}^2 dt
\end{equation}
Define 
\begin{equation}\label{eq:5.70}
H_T = \{u\in {\mathcal {D}}^{-N_0} \; : \; 
\| \langle g( P)\rangle^{-1}u\|^2_{\HH}+ \int_0^T \|A e^{  i t P } 
 u_{n} \|_{\HH_1}^2 dt< + \infty \}
\end{equation}
with its natural norm (the definition makes sense because of  the weak continuity property of $Ae^{itP}$). Due to the assumption~\eqref{eq:5.66} and the weak continuity property of $Ae^{itP}$, $H_T $ is a Hilbert space which is continuously embedded in ${\mathcal {D}}^{-N_{0}}$. The sequence $(u_{n})$ is bounded in $H$ and we can extract a subsequence converging weakly in $H$ to a limit $u$. 
Using the compactness of $(P+i)^{-1}$,  the operator $\bbbone_{ \langle
 P \rangle \leq \rho }$ is also compact on ${\mathcal {D}}^{-N_0}$. By passing to the 
limit we see that  $u$ satisfies 
\begin{equation}
C \| \bbbone_{ \langle
 P \rangle \leq \rho } u \|_{\HH}^2\geq 1
\end{equation} and 
\begin{equation}
0= \int_0^T \|A e^{  i t P } 
 u \|_{\HH_1}^2 dt
\end{equation}
The contradiction comes from the following:
\begin{lem}Denote by 
\begin{equation}
N=\{u \in H _T \; : \; 0= \int_0^T \|A e^{i t P } 
u \|_{\HH_1}^2 dt \}
\end{equation}
Then $N=\{0\}$.
\end{lem}
\begin{proof}
We first show that $N$ is invariant under the action of the operator $P$. 
Using that $Pe^{itP}u= i \partial _{t} e^{i t P } u$, the only thing to show is that 
if $u\in N$ then  $\| \langle g( P)\rangle^{-1}Pu\|_{\HH}$ is bounded.\par

We denote by $v(t) = e^{itP }u$ and apply~\eqref{eq5.new} with $T$ replaced by $T- \varepsilon_{0}$ to the sequence of functions 
\begin{equation}
v_{\varepsilon}= i\frac {v(t+ \varepsilon) - v(t)}\varepsilon
\end{equation}
we get for $0<\varepsilon< \varepsilon_{0}$ 
\begin{equation}\label{eq:5.new}
\| \langle g( P)\rangle^{-1}v_{\varepsilon}\rest_{t=0}\|_{\HH}\leq  C\| \bbbone_{ \langle
 P \rangle \leq \rho } v_{\varepsilon}\rest_{t=0} \|_{\HH}^2
\end{equation} 
and using that $v_{\varepsilon}\rest_{t=0}$ converges to $i\partial_{t}u\rest_{t=0}=Pu\rest_{t=0} $ in ${\mathcal {D}}^{-N_{0}-1}$, we obtain that the right hand side is bounded as $\varepsilon$ tends to $0$. Consequently, we can extract a subsequence $v_{\varepsilon}$ converging in $H_{T- \varepsilon_{0}}$. 
The limit is necessarily (due to the weak continuity property) $Pu$ which implies that $Pu\in N$.
To conclude, remark that $ \| \bbbone_{ \langle
 P \rangle \leq \rho } u_{n} \|_{\HH}^2$ is a norm on $N$ equivalent to the natural norm. Consequently $N$ is finite dimensional. The space $N$ is invariant by the operator $P$ which consequently has an eigenvector. But any eigenvector of $P$ in $N$ satisfies $Au=0$ and is equal to $0$ due to~\eqref{eq:5.1}. Consequently $N= \{0\}$.
\end{proof}
\section{Examples and applications}
\label{eaa}

In this section we present several applications of our method, giving, in particular
the proof of Theorems \ref{t:burq}, \ref{t:cdvp} and~\ref{t:bunim} stated in the introduction.

\subsection{Geometric control}\label{geomcont}
As in the introduction we consider 
$\Omega$, a smooth domain in ${\mathbb R}^d$,  $\Gamma\subset \partial \Omega$, and we fix $T>0$. 
For any $g\in L^2([0, T]\times  \Gamma)$, we denote by $u= S(g)$ the solution of the mixed problem \eqref{eq:1.1}.
The goal is to find conditions on 
$ \Gamma $ so that there exists 
a large class of functions $u_{0}$ which can be ``controlled'' by $g$, in the sense that
\begin{equation} 
\label{eq:6.control} 
u\rest_{t=T}=0 \,. 
\end{equation}
The basic result was obtained by Lebeau \cite{Le92} (see also~\cite{Li88} and~\cite{Zu88}). It
involves the natural concepts of the broken geodesic flow and of non-diffractive points (see
\cite{MeSj}, and also \cite{BuIMRN}):
\begin{thm}
\label{th1} Suppose that $\Gamma$ controls  $\Omega$ geometrically, 
that is 
\begin{equation}
\label{eq:gc}
\text{ $ \exists \; L_0 $ such that every trajectory of length $ L_0 $ meets $ \Gamma $ at a 
non-diffractive point,}
\end{equation}
where trajectories are with respect to the {\em broken geodesic flow}.
 Then for any $T>0$ and any $u_{1}\in H^1_{0}( \Omega)$ there exists 
$g\in L^2( [0,T] \times \Gamma)$ 
such that $S(g)\rest_{t>T}\equiv 0$.
\end{thm}
\begin{proof}
We first recall that as an application of Lions's H.U.M. method \cite{Li88} 
we see that  Theorem~\ref{th1} is equivalent to 
\begin{equation}\label{eq1.1}
\exists C>0; \quad \|u_{0}\|_{{H^1_{0}( \Omega)}}\leq C \|\partial_{n}(e^{it \Delta_{D}}u_{0})\rest_{[0,T]\times \Gamma}\|_{L^2( [0, T]\times \Gamma)}
\end{equation}
This follows from Theorem \ref{t:5.2} and the following resolvent estimate:
\begin{equation}
\label{eq:6.2}
\|R(z) f\|_{H^1( \Omega)} + \sqrt{ |z|}\|R(z) f\|_{L^2( \Omega)}\leq C  \|\partial _{n}R(z) f\|_{L^2( \Gamma)}+ 
 {C}  \| f\|_{L^2( \Omega)} \,,
\end{equation}
where $ R ( z ) = ( - \Delta_D - z )^{-1} $, with $ \Delta_D $, the Dirichlet Laplacian on $ \Omega $.
In fact, we can simply put $ A u = \partial_n u \rest_ \Gamma $ and $ \HH_1 = L^2 ( \Gamma ) $. 
To establish \eqref{eq:6.2} we can use the microlocal defect measures arguments as in \cite{BuIMRN}: we first prove ~\eqref{eq:6.2} for large $z$ and argue by contradiction. We obtain sequences $z_{n}\rightarrow + \infty$ and $u_{n}$ solution of 
\begin{gather}
\label{eq:6.2bis}
(-\Delta-z_{n}) u_{n}= f_{n}, \qquad \|u_{n}\|_{L^2( \Omega)}+\frac 1 {\sqrt {z_{n}}}\| \nabla _{x} u_{n}\|_{L^2( \Omega)}=1,
\\
\label{eq:6.2quart}
\|f_{n}\|_{L^2( \Omega)}
=o\left(\frac 1 {\sqrt{ z_{n}}}\right)
\\
\label{eq:6.2quint}
 \|\partial _{n} f\|_{L^2( \Gamma)} = o\left(\frac 1 {\sqrt{ z_{n}}}\right)
\end{gather}
Denote by $h_{n}= \sqrt{z_{n}}^{-1}$. Then, modulo the extraction of a
subsequence (see~\cite{GeLe93, BuIMRN}), 
there exists a positive Radon measure 
(a semi-classical defect measure) on $T^*{\mathbb R}^d$ such that, 
if $\underline{u_{n}}$ is the extension 
of $u_{n}$ by $0$ outside of $\overline{ \Omega}$, we have
\begin{enumerate}
\item For any $h$-pseudodifferential operator, $ A$,  on ${\mathbb R}^d$, we have
$$\langle \mu, \sigma_{0}(A)\rangle = \lim_{n\rightarrow + \infty} 
\left( A(x, h_{n}D_{x})\underline{u_{n}}, \underline{u_{n}}\right) _{L^2( {\mathbb R}^d}$$
\item The measure $\mu$ is supported in the semi-classical characteristic variety:
\begin{equation}\label{eq6.100}
\text{supp}( \mu) \subset T^*{\mathbb R}^d \cap \{ (x, \xi); x \in \overline{M},|\xi|^2= 1\end{equation}
\end{enumerate}
Furthermore (see~\cite{BuIMRN, BuLe01}),  using~\eqref{eq:6.2quart} we obtain that this measure is invariant along the generalized bicharacteristic flow. In the interior, this property is straightforward, whereas, near the boundary, it is more involved. In particular, we can show that the measure of the hyperbolic set (corresponding to transversal reflections) is equal to $ 0$.
This allows a definition of 
a bicharacteristic flow on the set~\eqref{eq6.100}, $\mu$ almost everywhere.  Due
to~\eqref{eq:6.2quint} the measure is equal to $0$ near any non diffractive point in $\Gamma$ (see~\cite{BuGe}); which, by~\eqref{eq:gc} implies that  the measure is identically null. Finaly the contradiction arises from the fact that according to~\eqref{eq:6.2bis} the measure has total mass $1$.\par
 The proof of~\eqref{eq:6.2} for $z\leq -1$ is straightforward using elliptic estimates and for $-1 \leq z \leq C$,~\eqref{eq:6.2} is obtained by a contradiction argument (and compactness) and the classical uniqueness theorem for second order elliptic operators (for this point we simply use that $\Gamma \neq \emptyset$).   
\end{proof} 

\subsection{Ikawa's black box}\label{ikba}
In the proof of Lebeau's theorem we did not use any ``black-box'' technology. As illustrated 
by Fig.\ref{fig:burq} we can employ it in 

\noindent
{\em Proof of Theorem \ref{t:burq}:}
As in the proof of Theorem \ref{th1} we use H.U.M. method and Theorem \ref{t:5.2} to 
reduce the argument to the following estimate:
\[
\|R(z) f\|_{H^1( \Omega)}+\sqrt{ |z|}\|R(z) f\|_{L^2( \Omega)}\leq C \log( |z|)\left(\|\partial _{n}R(z) f\|_{L^2( \Gamma)}+ \| f\|_{L^2( \Omega)}\right) \,,
\]
for $ \Im z \neq 0 $. This follows from Theorem \ref{t:5.1} and the following consequence of the
work of Ikawa~\cite[Theorem 2.1]{Ik}. Suppose that $ R_{\rm{bb}} ( k ) $ is the {\em outgoing}\footnote{The outgoing resolvent is the meromorphic continuation of $ (- \Delta - k^2 )^{-1} $ from $ \Im k > 0 $.}
resolvent for the
Dirichlet problem in the exterior of the union of convex obstacles satisfying
\begin{itemize}
\item $ ({\rm{convhull}}\; \Theta_j \cup \Theta_k  ) \cap \Theta_l = \emptyset \,, \ \ 
j \neq l \neq k \,.$
\item Denote by $\kappa$ the infimum of the principal curvatures of the boundaries of the obstacles $\Theta_{i}$, and $L$ the infimum of the distances between two obstacles. Then if $N>2$ we assume that $\kappa L >N$ (no assumption if $N=2$). 
\end{itemize}
Then there exist $ \alpha >0 $, $ C_0 $, and $ N_0 $ such that for $ \Im k > - \alpha $ we have
\[ \| \chi R_{\rm{bb}} ( k ) \chi \|_{L^2 \rightarrow L^2} \leq C_0 \langle k \rangle^{N_0} \,,
\ \ \chi \in \CI_c ( \RR^n ) \,.\]
An application of the maximum principle as in 
\cite[Lemma 2]{TZ98} and \cite[Lemma 4.10]{Bu02} (see also Lemma \ref{l:a3} below) gives a bound
\begin{equation}
\label{eq:rrr}
 \| \chi R_{\rm{bb}} ( k ) \chi \|_{L^2 \rightarrow L^2} \leq C_1 \frac{ \log \langle k
\rangle }{\langle k \rangle } \,,\end{equation}
and that gives the ``black-box'' assumption \eqref{eq:blackboxcontrol}  with $ G ( \lambda ) = \log \langle \lambda \rangle $ and $ A_{\rm{bb}} \equiv 0 $.
\stopthm

\subsection{Bunimovich stadium with the flat part as the black box}
\label{s:buni}
Our next control theoretical application is a new result about high frequency scarring in the case of the
Bunimovich 
stadium\footnote{which is perhaps the most celebrated example of a convex chaotic billiard}.
The same argument applies also in recent examples related to {\em quantum unique 
ergodicity} \cite{Do},\cite{Ze} where the flat part ``black box'' needs to be replaced by 
a flat torus. The result  which we use in the 
black box  (see Proposition \ref{prop:6.1} below) 
applies to that case as well. 

\begin{figure}[ht]
$$\ecriture{\includegraphics[width=12cm]{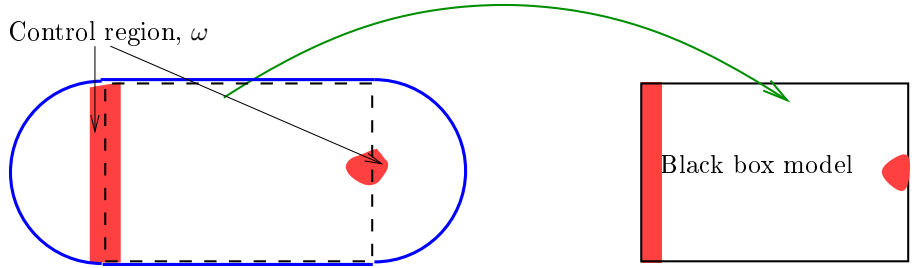}}
{\aat{1}{12}{Control region,  $\omega$}\aat{36}{5}{Black box model}
}$$
\caption{Control on the Bunimovich stadium}
\label{fig:bun}
\end{figure} 

\medskip
\noindent{\bf Theorem \ref{t:bunim}$'$.}{\em 
\ Consider $\Omega$ the Bunimovich stadium associated to a rectangle $R$, and 
$\omega\subset \Omega$ which controls a neighbourhood of 
$\Omega\setminus R$ geometrically. 
For any solution of the equation $(\Delta-z)v =f$, $ u \rest_{\partial \Omega} 
= 0 $ we then have 
\begin{equation}
\label{eq:6.11}\|v\|_{L^2( \Omega)}\leq C\left( \|f\|_{L^2( \Omega)}+ \|\bbbone_{\omega}v\|_{L^2( \omega)}\right)
\end{equation}}

We immediately deduce the following as a consequence of Theorem~\ref{t:5.2}: 
\begin{thm}
\label{t:bun}
Consider $\Omega$ the Bunimovich stadium associated to a rectangle $R$, and $\omega\subset \Omega$ which controls $\Omega\setminus R$
 geometrically. 
Then there exist $T>0$ and $C>0$ such that 
\begin{equation}
\label{eq:6.10}\|u_{0}\|^2_{{L^2( \Omega)}}\leq C \int _{0}^T \|\bbbone_{\omega}e^{it \Delta}u_{0}\|^2_{{L^2( \Omega)}}
\end{equation}
\end{thm}
In fact, by using a {\em temporal} black box, we could prove Theorem~\ref{t:bun} {\em for any} $T>0$.

We are going to deduce Theorem~\ref{t:bun}${}'$
 from the following result~\cite{Bu92} which is related to some earlier control results of Haraux~\cite{Ha} and  Jaffard~\cite{Ja}
\begin{prop}
\label{prop:6.1} Let $\Delta$ be  the Dirichlet 
Laplace operator on the rectangle $R= [0, 1]_{x} \times [0,a]_{y}$. 
Then for any open non-empty  $\omega\subset R$ of the form $
\omega= \omega_{x} \times [0,a]_{y}$ , there 
exists $C$ such that for any solutions of
\begin{equation}
(\Delta -z) u =f \ \text{ on $R$}\,, \ u \rest_{\partial R}=0 \,, \ z \in \RR \,, 
\end{equation}
we have
\begin{equation}
\label{eq:6.12}\|u\|^2_{{L^2(R)}}\leq C \left(\|f \|^2_{H^{-1} (
[0,1]_{x}; L^2([0,a]_{y})) } +
\|u\rest _{\omega}\|^2_{{L^2(\omega)}} \right)
\end{equation}
\end{prop}
\begin{proof}
We 
decompose $u,f$ in terms of 
the basis of $L^2([0,a])$ formed by the Dirichlet eigenfunctions
$e_{k}(y)=  { \sqrt {{2}/a}}\sin(2k\pi y/a)$,
\begin{equation}
u(x,y)= \sum_{k}e_{k}(y) u_{k}(x), \qquad f(x,y)= \sum_{k}e_{k}(y) f_{k}(x)
\end{equation}
we get for $u_{k}, f_{k}$ the equation
\begin{equation}\label{estres.1}
\left(\Delta_{x}-\left(z+\left({2k\pi}/{a}\right)^2\right)\right)u_{k}= f_{k},\qquad u_{k}(0)=u_k(1)=0
\end{equation}
Since $\omega_{x}$ controls geometrically $[0,1]$, a slight variant of~\eqref{eq:6.2} (or, in this simple case, a direct calculation) gives
\begin{equation}
\|u_{k}\|^2_{{L^2([0,1]_{x})}}\leq C \left(\|f_k \|^2_{H^{-1}([0,1]_{x})} +
\|u_k \rest _{\omega_{x}}\|^2_{{L^2(\omega)}}\right) 
\end{equation}
summing the squares on $k$ we get~\eqref{eq:6.12}.\footnote{We remark that as noted in \cite{Bu92} the proof applies to
any product manifold $M= M_{x}\times M_{y}$, and consequently Theorem 3' holds also for that geometry as a black box.}
\end{proof}

\medskip
\noindent
{\em Proof of Theorem \ref{t:bunim}$\,'$.}
Let us take $x,y$ as 
the coordinates on the stadium, so that $x$ is the longitudinal direction, 
$y$ the transversal direction, and the internal rectangle is 
$[0,1]_{x}\times [0,a]_{y}$. Let us then consider $u,f$ satisfying
 $(\Delta-z)u =f$, $ u = 0 $  on the boundary of
the stadium, and $\chi(x)\in C^\infty_{0}(0,1)$ equal to $1$ on $[\varepsilon, 1- \varepsilon]$. Then $\chi(x) u(x,y)$ is solution of
\begin{equation}
(\Delta-z)\chi u =\chi f + [\Delta, \chi] u\text{ in $R$}
\end{equation} with Dirichlet boundary conditions on $\partial R$. Applying Proposition~\ref{prop:6.1}, we get 
\begin{equation}
\|\chi u\|_{L^2( R)}\leq C\left( \|\chi f\|_{H^{-1}_{x}; L^2_{y}}+ \| u\rest _{\omega_{\varepsilon}}\|_{L^2( \omega_{\varepsilon})}\right)
\end{equation}
where $\omega_{\varepsilon}$ is a neighbourhood of the support of $\nabla \chi$.
Consequently we get for $V$ a neighbourhood of $\Omega\setminus R$, 
\begin{equation}\label{eqpresque}
\| u\|_{L^2( R)}\leq C\left( \| f\|_{L^2(R)}+ \| u\rest _{V}\|_{L^2(V)}\right)
\end{equation}
Finally, by standard propagation of semi-classical singularities as in Sect.\ref{geomcont}, we can replace in~\eqref{eqpresque} $V$ by $\omega$.  
\stopthm

\subsection{Semi-classical control with a prescribed loss.}
For completeness 
we present a natural 
class of 
examples in which $ G ( h ) $ in Theorems
\ref{t:3} and \ref{t:4} can essentially be a power of $ h $:
\[ G( h ) =  h^{-\alpha } \log (1/h) \,, \ \  \alpha = \frac{m-1}{m+1}\,, \ \  m=1,2,\cdots \,.  \]
For that consider the following set of Schr\"odinger operators on $ \RR^2 $:
\[ P_m ( h ) = - h^2 \Delta + x_1^2 - x_2^{2m} \,, \ \ m \in \NN  \,. \]
The Helffer-Sj\"ostrand theory of resonances \cite{HeSj} applies to this case 
(see also \cite[Sect.1]{SjDuke} where a discussion of a general polynomial is given).
In particular, for the meromorphically continued resolvent, $ R_m ( z, h ) = 
( P_m( h ) - z)^{-1} $, we have the following bound for the cut-off resolvent:
\begin{equation}
\label{eq:mb}  \| \chi R_m ( z , h ) \chi \| \leq C h^{-\frac{2m}{m+1}} \log(1/h)  \,.\end{equation} 
In fact, a separation of variables argument and the rescaling $ x = h^{\frac1{m+1}}y $ show
that the resonances are at the distance $ h^{\frac{2m}{m+1} } $ from the real axis. 
The same method shows that the resolvent is polynomially bounded in $ h^{-1} $ and
hence the interpolation argument we used before gives \eqref{eq:mb}.

>From $ P_m ( h ) $ we can construct a ``black box'' for an operator
 $ P ( h ) $ to which Theorems \ref{t:3} and \ref{t:4} will be 
applicable with $ G ( h ) = h^{-\frac{m-1}{m+1}} \log (1/h ) $.

\subsection{Closed hyperbolic orbits on manifolds.}
We will now discuss the case 
occuring when the black box contains a  hyperbolic orbit
in more detail, leading to the proof of Theorem \ref{t:cdvp}.

Thus suppose that the hypotheses of that theorem are satisfied.
It is well known that 
we can find a coordinate system
in a neighbourhood of $ \gamma $, $ U \simeq \SP^1 \times V $, 
$ V $ a neighbourhood of $ 0 $ in $ \RR^{n-1}$, in which $ \gamma $ 
is identified with $ \SP^1 $ and the metric is given by 
$$ g = d\theta^2 + \sum_{1\leq i,j \leq n-1} h_{ij} (x,\theta) dx_i 
dx_j \,, \ \  \theta \in \SP^1 \\, \ x \in V \,. $$
Since $ \gamma $ is hyperbolic we can assume that $ \SP^1 $ is the
only closed geodesic in $ U $. 

>From this local construction we now build a global scattering 
problem by extending $ g $ to a metric, $ g_{\rm{bb}} $, defined on 
$ \SP^1 \times \RR^{n-1} \simeq \SP^1_\theta \times \SP^{n-1}_\omega 
\times [0, \infty ) $. We choose $ g $ to be asymptotically Euclidean:
\[  g_{\rm{bb}} \sim dr^2 + r^2 d\theta^2 + r^2 g_{\SP^{n-1}} ( d\omega ) 
\,, \ \ r \rightarrow \infty \,,\]
and so that $ \gamma $ is the only closed geodesic of $ g_{\rm{bb}} $.

Because of the work of Ikawa \cite{Ik}, G\'erard \cite{Ge},
and of G\'erard-Sj\"ostrand \cite{GeSj},
it is expected that the  resolvent of the Laplacian of $ g_{\rm{bb}} $ can be
controlled using \eqref{eq:rrr},
as in Subsection \ref{ikba}. Since the two metrics agree in a
neighbourhood of the closed geodecics, we can use the scattering problem as
our ``black box'' and apply Theorem \ref{t:4} with $ A = ( 1 - \chi ) $. 
That would give Theorem \ref{t:cdvp} with $ (\log \lambda)^2 $ in place of 
$ \log \lambda $. To get the improved (and, thanks to an example in \cite{CdVP},
optimal) statement we need an improved estimate for the resolvent so that
Theorem \ref{t:4}$'$ can be applied:
\[ \| \chi R_{\rm{bb}} ( k ) \phi \|_{L^2 \rightarrow L^2} \leq C_1 \frac{ \sqrt{ \log \langle k
\rangle }}{\langle k \rangle } \,, \ \ \phi \in \CIc \,, \ 
\supp \phi \cap \gamma = \emptyset \,.\]

Since the needed results from scattering theory, although expected,
are not yet available\footnote{In \cite{Ik}
only convex obstacles in the Euclidean case are studied, while in \cite{GeSj} 
an analyticity assumption is made.} we take a simplified route and 
use a complex absorbing potential to construct a black box operator $ Q $ in
Theorem \ref{t:4}$'$\footnote{We remark however that the results 
of \cite{LG} and \cite{ChZw} would have
been sufficient for the case of hyperbolic geodesics on constant
negative curvature segments, if one takes the black box approach.}.
That is done in the Appendix with Theorem A furnishing
us with the needed estimates. Since we can use a neighbourhood of the
hyperbolic orbit of any Hamiltonian in phase space, we obtain a more general,
fully semi-classical variant of Theorem \ref{t:cdvp}:

\medskip
\noindent
{\bf Theorem \ref{t:cdvp}$'$.}\; {\em 
Suppose that  $ X $ is a compact $n$-manifold or $ \RR^n $, and $ P ( h ) \in \Psi_h^{m,0} ( X , 
\Omega_X^{\frac12} ) $ has
the principal symbol, $ p $,  satisfying: 
\begin{gather*}
 p^{-1} ([ - \epsilon , \epsilon ] ) \Subset T^* X \,, \ \text{ for some $ \epsilon > 0 $,}\\
p ( \rho ) = 0 \  \Longrightarrow \ dp ( \rho ) \neq 0 \,, \\ 
\exists \; C > 0 \ \ \  \langle \xi \rangle \geq C \ \Longrightarrow \  p \geq \langle \xi \rangle^m/C  \,,  
\end{gather*}
Let $ \gamma \subset p^{-1} ( 0 ) $ be closed hyperbolic orbit of the Hamilton flow
of $ p $, in the sense that all eigenvalues of the linearized Poincar\'e map are 
real and different from one. 

There exist constants $ C_0 $ and $ h_0 $, such that if 
$ u ( h ) \in L^2 ( X, \Omega_X^{\frac12} ) $ satisfies
\[  P( h ) u  = f  \,.\]
then for any $ A(h) \in \Psi_h^{0,0} ( X , \Omega_X^{\frac12} ) $, with 
its essential support, $ WF( A ) $, contained in a small neighbourhood of $ \gamma$, we have
\[ C_0\left( h^{-2} (\log(1/h))^2 \int_{X} |f|^2 + {\log(1/h)}\int_X | (I-A (h)) u  |^2\right) \geq \int_X | u  |^2 \,, 
\ \ h < h_0 \,.\]
}

\section*{Appendix}
\setcounter{section}{1}
\setcounter{equation}{0}
\setcounter{prop}{0}
\setcounter{lem}{0} 
\renewcommand{\thesection}{\Alph{section}}

In this appendix we will construct an operator $ Q $ appearing in Theorem \ref{t:4} for
a black box containing a hyperbolic orbit on a Riemannian manifold. Ideally, we 
would like $ Q $ to be the complex scaled Laplacian, $ -h^2 \Delta_\theta - z $ on an
asymptotically Euclidean manifold having one closed hyperbolic geodesic as its trapped
set. The results of \cite{Ge},\cite{GeSj} indicate that precise estimates of the type needed,
and in fact, the full understanding of resonances in logarithmic neighbourhoods of
the real axis, should be possible. Since we are dealing with the $ \CI $ case 
we will indicate here how the arguments of \cite{Ge} apply to this case. 

Let $ ( X, g ) $ be a scattering manifold satisfying the assumptions of \cite{WZ}.
In our application that means that near infinity $ X \simeq ( 0 , \epsilon ]_x \times
\SP^{n-2}_\omega \times \SP^{1}_\theta $, and the metric is $ g = dx^2/x^4 + 
g_{\SP^{n-2}}/x^2 + d\theta^2 / x^2  $, with infinity corresponding to $ x = 0 $. 
We assume that $ \gamma \subset X $ is the only closed geodesic on $ X $ and that
it is hyperbolic. 

Let $ a \in \CI ( X , [0,1] ) $ be equal to $ 0 $ in a neighbourhood of $ \gamma $ and
to $ 1 $, in a neighbourhood of infinity. We then put
\begin{equation}
\label{eq:a.q}
Q = Q ( z )  \stackrel{\rm{def}}{=} - h^2 \Delta_g - z - i h a \,,   \ \ 
z \in [1,2] + i [ - \epsilon , \epsilon ] \,.  
\end{equation}

The following result will allow applications of Theorem \ref{t:4}:

\medskip
\noindent
{\bf Theorem A.}{\em \; If $ Q ( z ) $ is given by \eqref{eq:a.q} and $ z \in I \Subset
(0, \infty) $, then for $ h < h_0 $, we have 
\begin{equation}
\label{t:a}
\| Q ( z )^{-1} \|_{ L^2 ( X ) \rightarrow L^2 ( X ) } 
\leq C \frac{ \log ( 1/ h ) }{h} \,.\end{equation}
If $ \phi \in \CIb ( X ) $ is supported away from $ \gamma $ then we also have
\begin{equation}
\label{t:a'}
\| Q ( z )^{-1} \phi \|_{ L^2 ( X ) \rightarrow L^2 ( X ) } 
\leq C \frac{ \sqrt{\log ( 1/ h )} }{h} \,.\end{equation}

}

To prove this theorem we will use the strategy of the proof of Theorem \ref{t:4} 
which means that it will be reduced to a local estimate near $ \gamma $.
We start with the well known version of Egorov's theorem. To state it 
we introduce an operator $ P \in \Psi^{m,0} ( X ) $ such that
\[ P = p^w( x , h D_x ; h ) + i h a ( x ) \,,  \ \ 
p ( x , \xi ;  h ) \in S^{m,0} ( T^* X ; \RR ) \,,
\ \ a \in \CIb ( T^*\RR ; \RR) \,. \]
We assume that the principal symbol of $ p $ satisfies $ p ( x , \xi ) 
\geq \langle \xi \rangle^m/C $ for $ |\xi| $ large enough. 
Then $ \exp ( - i t P / h ) $ is well defined and bounded on $ L^2 ( X ) $ 
either by the Hille-Yosida theorem or by a direct argument.
\begin{lem}
\label{l:a2}
Suppose that $ \Omega \subset \overline \Omega \Subset T^*X $, $ p \in S^{m,0} ( 
T^* X ) $ is real, and $ dp\rest_{ p^{-1}(0)} \neq 0 $ in $ \overline \Omega $. 
Suppose also that $ U \subset \Omega $ and that $ \exp ( t H_p ) U \subset \Omega $ for
$ 0 < t < T$. If $ p $ is the principal symbol of $ P \in \Psi^{m,0} $ and 
$ WF ( A ) $ is contained in $ U $, $ A \in \Psi^{0,0} ( T^* X ) $, $ \sigma_{m,0} ( A ) 
= a  $ then 
\begin{gather}
\label{eq:egor}
\begin{gathered}
 \exp (i t P/h ) A \exp ( - it P/h ) = \Op ( (\exp ( t H_p ) )^* a ) + E (t) \,, \\
 \| E (t) \|_{L^2 \rightarrow L^2 } \leq C_1 m ( A ) e^{ C_2 t  } h \,, \ 0 < t < T \,,
\end{gathered}
\end{gather}
where $ m ( A ) $ depends on a finite number of seminorms of the full symbol of $ A $, 
and $ C_1$, $ C_2 $ depend only on $ \Omega $ and $ p $.
\end{lem}
\noindent {\em Outline of the proof.} 
Using Proposition \ref{p:2.1} the result is obvious for 
$ U $ small enough and 
$ t $ such that  $ \bigcup_{ 0 \leq s\leq t } \exp ( s H_p ) U $
is contained in a sufficiently small neighbourhood of U. Since $ \Omega $ is precompact,
the size of $ U$ and $ t $ can be fixed uniformly in $ \Omega $. Assuming (as by a 
partitition of unity we may) that the $ U $ in the lemma is this small, we can divide
the interval $ [0,T] $ into subintervals of desired smallness. The errors estimates,
that is estimates on $ E( t ) $ in \eqref{eq:egor},
are multiplicative when switching from one interval to another and that gives the
exponential upper bound in $ t$.
\stopthm

\begin{figure}[ht]
$$\ecriture{\includegraphics[width=10cm]{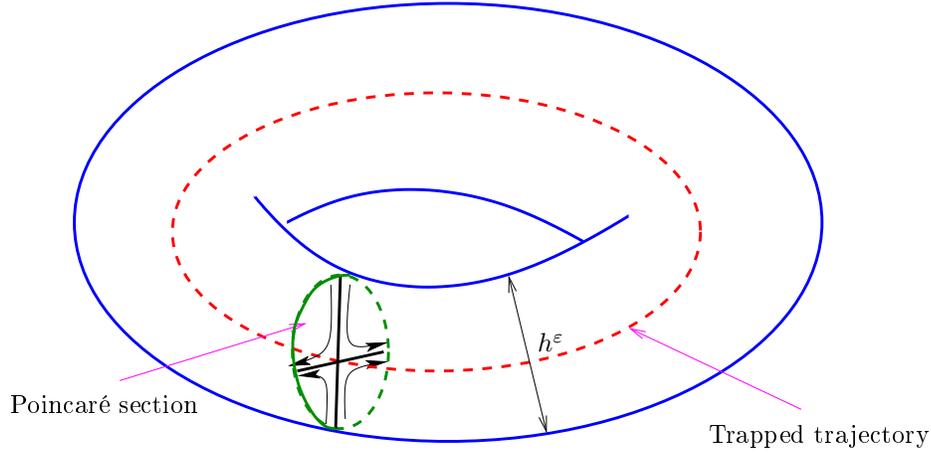}}
{\aat{42}{0}{Trapped trajectory}\aat{-3}{2}{Poincar\'e section}\aat{31}{6}{$h^\varepsilon$}
}$$
\caption{A hyperbolic trapped trajectory}
\label{fig:hyp}
\end{figure}

We can now show that we have control away from a small neighbourhood of
$ \gamma $. See Fig.\ref{fig:hyp} for an illustration of the hypotheses of the following
\begin{prop}
\label{p:a1}
Suppose that $ \epsilon $ is small, and let $ \psi_\epsilon \in S^0_\epsilon ( T^* X^\circ) 
\cap \CIc ( T^* X^\circ ) $ be a microlocal cut-off to an $ h^\epsilon$-neighbourhood of 
$ \pi^{-1} \gamma \cap \{ 1/2 \leq g ( x , \xi ) \leq 3 \} $, where $ g $ is
the metric.
Then, with  $ Q ( z ) $ as in \eqref{eq:a.q}, we have
\begin{equation}
\label{eq:al2}
Q ( z ) u = ( 1 - \psi_\epsilon ) f 
\ \Longrightarrow \ 
\| ( 1 - \psi_\epsilon ) u \|
 \leq  C \left( \frac{\log(1/h)}{h} \right) \| f \| + {\mathcal O} ( h^\infty ) \| u \| 
\,.
\end{equation}
If $ \epsilon = 0 $ then we have an improved estimate:
\begin{equation}
\label{eq:al22}
Q ( z ) u = ( 1 - \psi_0 ) f 
\ \Longrightarrow \ 
\| ( 1 - \psi_0 ) u \|
 \leq  C \frac{1}{h}  \| f \| + {\mathcal O} ( h^\infty ) \| u \| 
\,.
\end{equation}
\end{prop}
\begin{proof}
We will first prove \eqref{eq:al22} and then show how it implies \eqref{eq:al2}
using Lemma \ref{l:a2}. To see \eqref{eq:al22} we choose $ 
\tilde \psi_0 \in \CI $ so that $ ( 1 - \tilde \psi_0 ) ( 1 - \psi_0 ) = 
( 1 - \psi_0 ) $ and write
\[ \begin{split}
 h \int_X a |u|^2 & = 
 \Im \int_X Q ( z) u \overline u  = \int_X ( 1 - \psi_0 ) f \bar u \\
& \leq \| ( 1 - \psi_0 ) f \| \left( \| ( 1  - \tilde \psi_0 ) u \| 
+ {\mathcal O} 
(  h^\infty )   \| u \|\ \right) \,,
\end{split} \]
where we use the same symbols to denote the operator Weyl quantizing
the corresponding functions.
Lemma \ref{l:4} can be applied to $ Q ( z ) $ since 
both the imaginary term $ i a ( x ) h $ and $ z $ are lower order terms, and 
we can choose $  A u \stackrel{\rm{def}}{=} a ( x) u $. Hence
\[  \begin{split}
h \int_X a |u|^2 & \leq C \| ( 1 - \psi_0 ) f \| 
\left( \frac{1}{h} \| ( 1 - \psi_0 ) f \| + \|a u \| + {\mathcal O}(h^\infty ) 
\|u \| \right) \\
& \leq 2 \frac {C} { \varepsilon} h^{-1} \| ( 1 - \psi_0 ) f \|^2 + C \varepsilon h \| a u \|^2  + {\mathcal O}
(  h^\infty ) \| u \| 
\,, \end{split} \]
which proves \eqref{eq:al22}.

We now move to \eqref{eq:al2}. Let $ \varphi_\epsilon $ be a new microlocal cut-off
function localized to a annular neighbourhood, $ h^\epsilon < 
d ( \bullet, \gamma ) < h^{\epsilon/2} $. Splitting it into incoming and outgoing
parts with respect to the flow,  we can, by forward and retarded propagation respectively,
move it by $ \exp ( - i t Q ( h ) /h ) $, $ |t| \simeq \epsilon \log( 1 / h ) $
into a fixed size set, a finite distance from $ \gamma $ and away 
from the support of $ a $. The last condition guarantees that the propagator
is microlocally unitary. We can then apply \eqref{eq:al22}. We can continue
by a dyadic decomposition argument, with the number of terms proportional to 
$ \log ( 1/h) $.
\end{proof}

With the help of the above result  we have essentially reduced the proof of Theorem A to 
the proof of the following

\begin{prop}
\label{p:a2}
With the notation of Proposition \ref{p:a1} there exist $ c_0 $, $ h_0 $, and $ N_0 $ such that 
we have 
\begin{equation}
\label{eq:pa2}
 Q ( z ) u =  \psi_\epsilon  f 
\ \Longrightarrow \ 
\|  \psi_\epsilon  u \|
 \leq  C {h}^{-N_0}  \| f \| + {\mathcal O} ( h^\infty ) \| u \| 
\, \end{equation}
if $ z \in [1,2] + i(-c_0 h , +\infty) $ and $ h < h_0 $.
\end{prop}

\noindent {\em Outline of the proof.}
Using \cite[Proposition 5.1]{SjZw02} we can reduce the proof of \eqref{eq:pa2} 
to an estimate for an operator involving the quantum monodromy operator, $ M ( z ) $
(see \cite[Sect.4]{SjZw02}, and, for a brief introduction, \cite[Sect.2, Appendix]{ISZ}):
\begin{equation}
\label{eq:what}
 \| \psi^\sharp_\epsilon ( I - M ( z) )^{-1} \psi^\sharp_\epsilon \|_{ L^2( \RR^{n-1})
\rightarrow L^2 ( \RR^{n-1} ) } = {\mathcal O} ( h^{-N_0+1} ) \,, \  \
 z \in [1,2] + i(-c_0 h , c_0 h ) \,,\end{equation}
where $  \psi^\sharp_\epsilon $ is a microlocal cut-off to 
an $h^\epsilon$-neighbourhood of $ ( 0 , 0 ) \in T^* \RR^{n-1} $, induced by $ \psi_\epsilon $
after the identification with the Poincar\'e section (see Fig.\ref{fig:hyp}), and the
inverse of $  I - M ( z )  $ exists on a Hilbert space
 $ H^p_h \subset L^2 $, and such that 
\[  \psi_\epsilon^\sharp = {\mathcal 
O} ( h^{-N(p)} ) \; : \; L^2 \rightarrow H^p_h \,, \ \ 
 \psi_\epsilon^\sharp = {\mathcal 
O} ( h^{-N(p)} ) \; : \;  H^p_h  \rightarrow L^2 \,. \] 
The operator $ M ( z ) $ is of the same form as the operator $ e^{i z \ell ( \gamma ) /h } 
M_1 $, $ \ell(\gamma ) $ the length of $ \gamma $,  of 
\cite[Th\'eor\`eme 2.6]{Ge}. We need a translation from the large parameter 
setting of \cite{Ge} to the semiclassical setting: 
$ \lambda = z/h $, $ h = 1/\lambda_1$, $ \lambda_1 = \Re \lambda $. The spaces $ H^p_h$
are defined in \cite[\S 3.1, \S 4.2]{Ge} and the estimates for the
Grushin problem for $ I -  M ( z ) $ are obtained in \cite[Th\'eor\`eme 4.11]{Ge} (the
variable $ z $ appearing there is $ \exp ( i z \ell(\gamma) / h ) $ in our notation).
Since in \eqref{eq:what} we only need the invertibility of the resulting effective 
Hamiltonian for $ z \in [1,2] + i(-c_0 h , c_0 h ) $,
G\'erard's analysis gives that and much more.
\stopthm

To prove Theorem A we need the following lemma which, for possible future use,
we state in a slightly excessive generality:
\begin{lem}
\label{l:a3}
Suppose that $ A $ and $ B $ are 
bounded self-adjoint operators on a Hilbert space $ \HH $, 
\[ A^2 = A \,, \ \ B A = A B = A \,,\]
and $ F( z) $ is a family of bounded operators satisfying 
\begin{gather}
\label{eq:F}
\begin{gathered}
 F( z )^* = F( \bar z ) \,, \ \ \partial_z F\rest_\RR \; 
\geq c Id \,, \ c > 0 \,, \\
\text{ $ B F ( z )^{-1} B $ is holomorphic in $ [-\epsilon , \epsilon ] + i [-\delta , \delta ] 
 $,  \ $ \frac{\delta}{\epsilon} \ll 1/\log M $},  \\
\| B F ( z )^{-1} B \| \leq M \,, \ \  \| A F ( z )^{-1} A \| \leq 1 \,.
\end{gathered}
\end{gather}
Then for $ |z | < \epsilon /2 $, $ \Im z = 0 $ we have 
\begin{equation}
\label{eq:gain}
\| B F ( z )^{-1} B \| \leq C  \frac{\log M}\delta  \,, \ \
\| B F ( z )^{-1} A \| \leq C \sqrt{ \frac{\log M}\delta  } \,. 
\end{equation}
\end{lem}
\begin{proof}
The first part of \eqref{eq:gain} works exactly as in \cite[Lemma 2]{TZ98} and 
\cite[Lemma 4.2]{Bu02}. To see the improved version we 
 start by observing that the conditions on $ F $ and $ A $ imply that for
$ \Im z > 0 $, small, 
\[ \Im z  \| u \|^2 \leq C \Im \langle F( z ) u , u \rangle \,.\]
If now $ F (z) u = A f $, then by the assumptions on $ F $, $ \| A u \|
\leq \| A f \| $, and consequently,
\[ \| B u \|^2 \leq C \| u \|^2 \leq \frac{1}{ \Im z } \langle A f , A u \rangle 
\leq \frac{1}{\Im z } \| A f \|^2 \,,\]
Here we used the facts that $ A^2 = A = A^* $.
Since $ u = F ( z)^{-1} A f $,
this, and the fact that $ B A = A $, give 
\[ \begin{split}
& \| B F( z )^{-1} A \| \leq \frac{C }{\sqrt{\Im z }} \,, \ \ \Im z > 0 \\
& \| B F ( z)^{-1} A \| \leq C \| B F ( z) ^{-1} B \| \leq C M \,,
\end{split}\]
Interpolating as before gives \eqref{eq:gain}.
\end{proof}

\medskip
\noindent
{\em Proof of Theorem A}. 
We first combine Propositions~\ref{p:a1} and~\ref{p:a2} to estimate $(1-\psi_{\varepsilon}) Q^{-1} (1-\psi_{\varepsilon})$ and $\psi_{\varepsilon} Q^{-1} \psi_{\varepsilon}$ by $h^{-N}$. Then, since 
\[Q(1-\psi_{\varepsilon})Q^{-1}\psi_{\varepsilon}f=-[Q, \psi_{\varepsilon}]Q^{-1}\psi_{\varepsilon}f + (1-\psi_{\varepsilon})\psi_{\varepsilon}f\]
by using these estimates (for a different function $\psi$) we get an estimate of the same type for $(1-\psi_{\varepsilon}) Q^{-1} \psi_{\varepsilon}$ and consequently for $Q^{-1}$. Finally we combine this latter estimate and~\eqref{eq:al22} with 
Lemma \ref{l:a3} 
applied to the family of operators $ w \mapsto F ( w ) = (i/h)Q( z_0 + hw ) $ and 
$  A  =  \bbbone_{\supp \phi} $,   $ B =  1 $. We can take $ \delta $ independent of $ h$ and $ \epsilon = 1/(Ch) $ so that the assumption $ \delta/\epsilon 
\ll 1/\log(1/h^N) $ is easily satisfied.
\stopthm

\def\cprime{$'$}

\end{document}